\documentclass[11pt]{amsart}
\usepackage{amsfonts}
\usepackage{amsmath,amsthm,amsfonts,amssymb}
\newtheorem{lemma}{Lemma}
\newtheorem{theorem}{Theorem}
\newtheorem{corollary}{Corollary}
\newtheorem*{thmA}{Theorem A}
\newtheorem*{thmB}{Theorem B}
\newtheorem*{thmC}{Theorem C}

\def\bl{\begin{lemma}}
\def\bt{\begin{theorem}}
\def\el{\end{lemma}}
\def\et{\end{theorem}}
\def\bp{\begin{proof}}
\def\ep{\end{proof}}
\def\bc{\begin{corollary}}
\def\ec{\end{corollary}}
\def\iy{\infty}

\def\mc{\mathcal}
\def\mb{\mathbb}

\def\l{\lambda}

\def\o{\omega}
\def\O{\Omega}
\def\a{\alpha}

\def\p{\partial}

\def\s{\sigma}

\def\r{\rho}

\def\-{\setminus}
\def\s{\sigma}

\def\vp{\varphi}
\def\ov{\overline}
\def\lt{\left}
\def\rt{\right}
\def\+{\bigcup}
\def\.{\bigcap}

\def\ll{\langle}
\def\rl{\rangle}

\title[A Schwarz lemma for
harmonic mappings]{A Schwarz lemma for harmonic mappings between the
unit balls in real Euclidean spaces}
\thanks{Research supported by the National Natural Science Foundation of China
(No. 11201199) and by the Scientific Research Foundation of Jinling
Institute of Technology (No. Jit-b-201221).}

\author {Shaoyu Dai and Yifei Pan}

\address{Department of General Study Program, Jinling Institute of
Technology, Nanjing 211169, China}
\address{\it E-mail address: dymdsy@163.com}

\address{Department of Mathematical Sciences, Indiana University -
Purdue University Fort Wayne, Fort Wayne, IN 46805-1499, USA}
\address{\it E-mail address: pan@ipfw.edu}

\begin{document}

\numberwithin{equation}{section}

\begin{abstract}
In this paper we prove a Schwarz lemma for harmonic mappings between
the unit balls in real Euclidean spaces. Roughly speaking, our
result says that under a harmonic mapping between the unit balls in
real Euclidean spaces, the image of a smaller ball centered at
origin can be controlled. This extends the related result proved by
Chen in complex plane.
\end{abstract}

\maketitle



\section{Introduction}

Let $n$ be a positive integer greater than 1. $\mathbb{R}^n$ is the
real space of dimension $n$. For
$x=(x_1,\cdots,x_n)\in\mathbb{R}^n$, let
$|x|=(|x_1|^2+\cdots+|x_n|^2)^{1/2}$. Let
$\mathbb{B}^n=\{x\in\mathbb{R}^n: |x|<1\}$ be the unit ball of
$\mathbb{R}^n$. The unit sphere, the boundary of $\mathbb{B}^n$ is
denoted by $S$; normalized surface-area measure on $S$ is denoted by
$\s$ (so that $\s(S)=1$). Let $S^+$ denote the northern hemisphere
$\{x=(x_1,\cdots,x_n)\in S: x_n>0\}$ and let $S^-$ denote the
southern hemisphere $\{x=(x_1,\cdots,x_n)\in S: x_n<0\}$.
$N=(0,\cdots,0,1)$ denotes the north pole of $S$.
$B^n_r=\{x\in\mathbb{R}^n: |x|<r\}$ is the open ball centered at
origin of radius $r$; its closure is the closed ball $\ov {B^n_r}$.

Let $m$ be a positive integer with $m\geq1$. A mapping
$F=(F_1,\cdots,F_m,F_{m+1})$ from $\mathbb{B}^n$ into
$\mathbb{B}^{m+1}$ is harmonic on $\mathbb{B}^n$ if and only if for
$k=1,\cdots,m,m+1$, $F_k$ is twice continuously differentiable and
$\Delta F_k\equiv0$, where $\Delta=D_1^2+\cdots+D_n^2$ and $D_j^2$
denotes the second partial derivative with respect to the $j^{th}$
coordinate variable $x_j$. By $\O_{n,m+1}$, we denote the class of
all harmonic mappings $F$ from $\mb{B}^n$ into $\mb{B}^{m+1}$.

Let $\mathfrak{B}^n$ be the unit ball in the complex space
$\mathbb{C}^n$. Denote the ball $\{z\in\mathbb{C}^n: |z|<r\}$ by
$\mathfrak{B}^n_r$; its closure is the closed ball $\ov
{\mathfrak{B}^n_r}$. For a holomorphic mapping $f$ from
$\mathfrak{B}^n$ into $\mathfrak{B}^m$, the classical Schwarz lemma
\cite{Rudin} says that if $f(0)=0$, then
\begin{equation}\label{hljSlemma1}
|f(z)|\leq|z|
\end{equation}
holds for $z\in\mathfrak{B}^n$. For $0<r<1$, \eqref{hljSlemma1} may
be written in the following form:
\begin{equation*}
f(\ov {\mathfrak{B}^n_r})\subset\ov {\mathfrak{B}^m_r}.
\end{equation*}
So the classical Schwarz lemma can be regarded as considering the
region of $f(\ov {\mathfrak{B}^n_r})$. If $f(0)\neq0$, then what the
region of $f(\ov {\mathfrak{B}^n_r})$ is. It seems that there is not
much of research in the literature. However, the same problem also
exists in harmonic mappings. The work in the following by Chen
\cite{CHH} seems to be the first result of this kind of study for
harmonic mappings in the complex plane.

Let $\mathbb{D}$ be the unit disk in the complex plane $\mathbb{C}$.
Denote the disk $\{z\in\mathbb{C}: |z|<r\}$ by $D_r$; its closure is
the closed disk $\ov D_r$. For $0<r<1$ and $0\leq\r<1$, Chen
\cite{CHH} constructed a closed domain $E_{r,\r}$ and proved that
\begin{thmA}
Let $0\leq\r<1$, $\a\in\mathbb{R}$ and $0<r<1$ be given. For every
complex-valued harmonic function $F$ on $\mathbb{D}$ such that
$F(\mathbb{D})\subset\mathbb{D}$, if $F(0)=\r e^{i\a}$, then
\begin{equation}\label{hljChen}
F(\ov D_r)\subset e^{i\a}E_{r,\r},\
\end{equation}
which is sharp.
\end{thmA}

Note that the function $F$ in the above theorem can be seen as
$F\in\O_{2,2}$. So \eqref{hljChen} can be regarded as considering
the region of $F(\ov {B_r^2})$ when $F\in\O_{2,2}$ regardless of
$F(0)=0$ or $F(0)\neq0$. In \cite{CHH}, the most important theorem
for the proof of Theorem A is the theorem as follow, which is the
motivation for our study of the extremal mapping. The mappings
$U_{a,b,r}$ and $F_{a,b,r}$ in the following theorem are defined in
\cite{CHH}.
\begin{thmB} Let $F=U+iV$ be a harmonic mapping such that
$F(\mathbb{D})\subset\mathbb{D}$ and $F(0)=a+bi$. Then for $0<r<1$
and $0\leq\theta\leq2\pi$,
$$U(re^{i\theta})\leq U_{a,b,r}(ri)$$
with equality at some point $re^{i\theta}$ if and only if
$F(z)=F_{a,b,r}(e^{i(\pi/2-\theta)}z)$. Furthermore,
$U(z)<U_{a,b,r}(ri)$ for $|z|<r$. \end{thmB}

A classical Schwarz lemma for complex-valued harmonic function on
$\mathbb{B}^n$ \cite{ABW} says that
\begin{thmC}
Suppose that $F$ is a complex-valued harmonic function on
$\mathbb{B}^n$, $|F|<1$ on $\mathbb{B}^n$, and $F(0)=0$. Then
\begin{equation}\label{ybds}
|F(x)|\leq U(|x|N)
\end{equation}
holds for every $x\in\mathbb{B}^n$, where $U$ is the Poisson
integral of the function that equals 1 on $S^+$ and -1 on $S^-$.
Equality holds for some nonzero $x\in\mathbb{B}^n$ if and only if
$F=\l(U\circ A)$ where $\l$ is a complex constant of modulus $1$ and
$A$ is an orthogonal transformation.
\end{thmC}

Especially, when $n=2$ in the above theorem, it is known
\cite{Heinz} that
$$|F(x)|\leq \frac{4}{\pi}\arctan|x|$$
holds for every $x\in\mathbb{B}^2$.

From Theorem C, for $0<r<1$,
\eqref{ybds} may be written in the following form:
\begin{equation}\label{hljHSlemma}
F(\ov {B_r^n})\subset\ov D_{U(rN)},
\end{equation}
where $\ov D_{U(rN)}=\{z\in\mathbb{C}: |z|\leq U(rN)\}$.

Note that the function $F$ in the above Theorem C can be seen as
$F\in\O_{n,2}$. So \eqref{hljHSlemma} can be regarded as considering
the region of $F(\ov {B_r^n})$ when $F\in\O_{n,2}$ with $F(0)=0$. It
is natural to consider that if $F\in\O_{n,2}$ with $F(0)\neq0$, then
what the region of $F(\ov {B_r^n})$ is. Furthermore, we want to know
that for the general $F\in\O_{n,m+1}$, what the estimate
corresponding to \eqref{hljHSlemma} is when $F(0)=0$ or $F(0)\neq0$.
This problem will be resolved in this paper. When $F(0)\neq0$, this
problem is serious because the composition $f\circ F$ of a
$\mbox{m$\ddot{o}$bius}$ transformation $f$ and a harmonic mapping
$F$ does not need to be harmonic.

In this paper, inspired by the method of the proof of Theorem B in
\cite{CHH}, we obtain the following Theorem \ref{e}, which is very
important in this paper. \eqref{zy1} is the estimate corresponding
to \eqref{ybds} without the assumption $F(0)=0$. Especially, when
$F(0)=0$, we have Corollary \ref{especial}, which is coincident with
Theorem C when $m+1=2$. Note that in the following theorem,
$F_{(a,b)Q_e,r}$ is defined as \eqref{Fabrdy}.

\bt\label{e} Let $F(x)$ be a harmonic mapping such that
$F(\mathbb{B}^n)\subset\mathbb{B}^{m+1}$ and $F(0)=(a,b)$, where
$a\in\mathbb{R}^{m}$ and $b\in\mathbb{R}$. Let $e$ be a unit vector
in $\mathbb{R}^{m+1}$, $e_0=(1,0,\cdots,0)\in\mathbb{R}^{m+1}$ and
$Q_e$ be an orthogonal matrix such that $eQ_e=e_0$. Then, for
$0<r<1$ and $\o\in S$,
\begin{equation}\label{zy1}
\ll F(r\o),e\rl\le \ll F_{(a,b)Q_e,r}(rN),e_0\rl
\end{equation}
with equality at some point $r\o$ if and only if
$F(x)=F_{(a,b)Q_e,r}(xA)Q_e^{-1}$, where $A$ is an orthogonal matrix
such that $\o A=N$ and $Q_e^{-1}$ is the inverse matrix of $Q_e$.
Furthermore, $\ll F(x),e\rl<\ll F_{(a,b)Q_e,r}(rN),e_0\rl$ for
$|x|<r$. \et

\bc\label{especial} Let $F(x)$ be a harmonic mapping such that
$F(\mathbb{B}^n)\subset\mathbb{B}^{m+1}$ and $F(0)=0$. Then
\begin{equation*}
|F(x)|\leq U(|x|N)
\end{equation*}
for every $x\in\mathbb{B}^n$, where $U$ is the Poisson integral of
the function that equals 1 on $S^+$ and -1 on $S^-$. Equality holds
for some nonzero $x_0\in\mathbb{B}^n$ if and only if $F(x)=U(xA)e$,
where $A$ is an orthogonal matrix such that $x_0A=|x_0|N$, $e$ is a
unit vector in $\mathbb{R}^{m+1}$. \ec

From Theorem \ref{e}, we deduce the following theorem, which is
called a harmonic Schwarz lemma for $F\in\O_{n,m+1}$ and which
resolves the problem we want to know above. Theorem \ref{egj}
extends Theorem A and is coincident with Theorem A when $n=m+1=2$.
Note that in the following theorem, $F_{(a,b)Q_e,r}$ is defined as
\eqref{Fabrdy}.

\bt\label{egj} Let $F(x)$ be a harmonic mapping such that
$F(\mathbb{B}^n)\subset\mathbb{B}^{m+1}$ and $F(0)=(a,b)$, where
$a\in\mathbb{R}^{m}$ and $b\in\mathbb{R}$. Let $0<r<1$. Then
\begin{equation}\label{zhjg}
F(\ov {B_r^n})\subset E_{r,(a,b)},
\end{equation}
where
$$E_{r,(a,b)}=\bigcap_{e\in\mathbb{R}^{m+1},|e|=1}R_e,$$
$$R_e=\{x\in\mathbb{R}^{m+1}: \ll x,e\rl\le \ll F_{(a,b)Q_e,r}(rN),e_0\rl\},$$
$e_0=(1,0,\cdots,0)\in\mathbb{R}^{m+1}$ and $Q_e$ be an orthogonal
matrix such that $eQ_e=e_0$. \et

Note that $E_{r,(a,b)}$ in Theorem \ref{egj} is a region enveloped
by all the hyperplanes
$$P_e=\{x\in\mathbb{R}^{m+1}: \ll x,e\rl=\ll
F_{(a,b)Q_e,r}(rN),e_0\rl\},$$ which is the boundary of $R_e$. By
Theorem \ref{e}, it is obviously that the region $E_{r,(a,b)}$ is
sharp. This means that under $F\in\O_{n,m+1}$, the image of a small
ball centered at origin of radius $r$ can be controlled.

In section \ref{lemmas}, we will give two main lemmas. The proofs of
the lemmas will be given in section \ref{The proofs}. In section
\ref{results}, the main results of this paper and the proofs will be
given.

\section{The main lemmas}\label{lemmas}

In this section, we will introduce two main lemmas, which are
important for the proof of Theorem \ref{theoremuabr} and which
extend the related lemmas proved by Chen in \cite{CHH}. Lemma
\ref{lemmaRI} constructs a bijection $(R,I)$ from
$\mathbb{R}^m\times\mathbb{R}^+$ onto the upper half ball
$\{(a,b):a\in\mathbb{R}^m,b\in\mathbb{R},|a|^2+b^2<1,\ b>0\}$, which
will be used to construct $u_{a,b,r}$ in Theorem \ref{theoremuabr}
for the case that $b>0$. Lemma \ref{lemmaR} constructs a bijection
$\mc R$ from $\mathbb{R}^m$ onto the ball $\{a: a\in\mathbb{R}^m,
|a|<1\}$, which will be used to construct $u_{a,b,r}$ in Theorem
\ref{theoremuabr} for the case that $b=0$. Now we give the two main
lemmas. The proofs of Lemma \ref{lemmaRI} and Lemma \ref{lemmaR}
will be given in section \ref{The proofs}.

For $0<r<1$, $\mu>0$, $\l\in\mathbb{R}^m$, and
$l=(1,0,\cdots,0)\in\mathbb{R}^m$, define
\begin{equation}\label{Ady}
A_{r,\l,\mu}(\omega) =\frac1\mu\lt(\frac{1}{|rN-\o|^n}l-\l\rt),\ \ \
\ \ \ \o\in S,
\end{equation}
and
\begin{equation}\label{RIdy}
R(r,\l,\mu)=\int_S\frac{A_{r,\l,\mu}(\o)}
{\sqrt{1+|A_{r,\l,\mu}(\o)|^2}}\,d\s, \quad
I(r,\l,\mu)=\int_S\frac{1} {\sqrt{1+|A_{r,\l,\mu}(\o)|^2}}\,d\s.
\end{equation}
The idea of the conformation of $A_{r,\l,\mu}(\omega)$,
$R(r,\l,\mu)$ and $I(r,\l,\mu)$ originates from \eqref{u0dy} and
\eqref{u4}.

\bl\label{lemmaRI} Let $0<r<1$ be fixed. Then, there exist a unique
pair of continuous mappings $\l=\l(r,a,b)\in\mathbb{R}^m$ and
$\mu=\mu(r,a,b)>0$, defined on the upper half ball
$\{(a,b):a\in\mathbb{R}^m,b\in\mathbb{R},|a|^2+b^2<1,\ b>0\}$, such
that $R(r,\l(r,a,b),\mu(r,a,b))=a$ and $I(r,\l(r,a,b),\mu(r,a,b))=b$
for any point $(a,b)$ in the half ball. \el

For $0<r<1$, $\l\in\mathbb{R}^m$, and
$l=(1,0,\cdots,0)\in\mathbb{R}^m$, define
\begin{equation}\label{A}
\mc A_{r,\l}(\omega)=\frac{1}{|rN-\o|^n}l-\l,\ \ \ \ \ \ \o\in S,
\end{equation}
and
\begin{equation}\label{R}
\mc R(r,\l)=\int_S\frac{\mc A_{r,\l}(\o)} {|\mc A_{r,\l}(\o)|}\,d\s.
\end{equation}
The idea of the conformation of $\mc A_{r,\l}(\omega)$ and $\mc
R(r,\l)$ originates from \eqref{b=0u2}. Note that $\mc R(r,\l)$ is
well defined, since $|\mc A_{r,\l}(\o)|\neq0$ except for a zero
measure set of $\o$ at most.

\bl\label{lemmaR} Let $0<r<1$ be fixed. Then, there exist a unique
continuous mapping $\l=\l(r,a)\in\mathbb{R}^m$, defined on
$\{a:a\in\mathbb{R}^m,|a|<1\}$, such that $\mc R(r,\l(r,a))=a$ for
any point $a$. \el

\section{The main results}\label{results}

Let $a\in\mathbb{R}^m$, $b\in\mathbb{R}$ and $0\le b<1$,
$|a|^2+b^2<1$. Let ${\mc U}_{a,b}$ denote the class of mappings
$u\in (L^\infty(S))^m$ satisfying the following conditions:
\begin{equation} \label{uyq}
\|u\|_{\infty}\le 1,\quad \int_S u(\o)d\s=a,\quad
\int_S\sqrt{1-|u(\o)|^2}d\s\ge b.
\end{equation}
Every function $u\in (L^\infty(S))^m$ defines a harmonic mapping
$$U(x)=\int_S\frac{1-|x|^2}{|x-\o|^n}
u(\o)d\s\quad\mbox{for}\quad x\in\mathbb{B}^n.$$ Let $0<r<1$,
$l=(1,0,\cdots,0)\in\mathbb{R}^m$ and define a functional $L_r$ on
$(L^\infty(S))^m$ by
\begin{equation} \label{Ldy}
L_r(u)=\ll U(rN),l\rl =\int_S\frac{1-r^2}{|rN-\o|^n} \ll u(\o),l\rl
d\s.
\end{equation}

Obviously, ${\mc U}_{a,b}$ is a closed set, and $L_r$ is a
continuous functional on ${\mc U}_{a,b}$. Then there exists a
extremal mapping such that $L_r$ attains its maximum on ${\mc
U}_{a,b}$ at the extremal mapping. We will claim in the following
theorem that the extremal mapping is unique. In the proof of the
following theorem, we will construct a mapping $u_0$ first and then
prove that $u_0$ is the unique extremal mapping, which will be
denoted by $u_{a,b,r}$.

\bt\label{theoremuabr} For any $a$, $b$ and $r$ satisfying the above
conditions, there exists a unique extremal mapping $u_{a,b,r}\in{\mc
U}_{a,b}$ such that $L_r$ attains its maximum on ${\mc U}_{a,b}$ at
$u_{a,b,r}$. \et

For the proof of Theorem \ref{theoremuabr}, we need Lemma 1, Lemma 2
and the lemma in the following.

\bl\label{lemmaxybds} Let $x,y\in\mathbb{R}^m$, $|x|\leq1$ and
$|y|<1$. Then
\begin{equation}\label{xy1}
\sqrt{1-|y|^2}-\sqrt{1-|x|^2}=\frac{\ll x-y,y\rl}{\sqrt{1-|y|^2}}
+\frac{|x-y|^2(1-|\tilde{y}|^2)+|\ll
x-y,\tilde{y}\rl|^2}{2(1-|\tilde{y}|^2)^{3/2}}
\end{equation}
holds, where $\tilde{y}=y+\zeta(x-y)$, $0<\zeta<1$. \el

\bp Let $x=(x_1,\cdots,x_m)$, $y=(y_1,\cdots,y_m)$ and
$g(x)=\sqrt{1-|x|^2}$. For $j=1,\cdots,m$ and $k=1,\cdots,m$, denote
$\frac{\p g(x)}{\p x_j}$ by $g_j(x)$ and $\frac{\p^2 g(x)}{\p x_j\p
x_k}$ by $g_{jk}(x)$. Then
$$g_j(x)=-\frac{x_j}{\sqrt{1-|x|^2}},$$
$$\frac{\p^2 g(x)}{\p x_j\p
x_k}=-\frac{\delta_{jk}(1-|x|^2)+x_jx_k}{(1-|x|^2)^{3/2}},$$ where
\begin{equation*}
\delta_{jk}=
\begin{cases}
1, &j=k;\\
0, &j\neq k.
\end{cases}
\end{equation*}
Let $\vp(t)=g(y+t(x-y))$. By Taylor formula, we have
\begin{equation} \label{xy2}
\vp(1)-\vp(0)=\vp'(0)+\frac{1}{2}\vp''(\zeta)\ \ \ \ \ \ \
(0<\zeta<1).
\end{equation}
Note that
$$\vp(0)=g(y)=\sqrt{1-|y|^2},\ \ \ \ \ \
\vp(1)=g(x)=\sqrt{1-|x|^2},$$
$$\vp'(t)=\sum_{j=1}^mg_j(y+t(x-y))\cdot(x_j-y_j),$$
$$\vp'(0)=\sum_{j=1}^mg_j(y)\cdot(x_j-y_j)=-\frac{\ll
x-y,y\rl}{\sqrt{1-|y|^2}},$$
$$\vp''(t)=\sum_{j,k=1}^mg_{jk}(y+t(x-y))\cdot(x_j-y_j)(x_k-y_k),$$
\begin{eqnarray*}
&&\vp''(\zeta)=\sum_{j,k=1}^mg_{jk}(\tilde{y})\cdot(x_j-y_j)(x_k-y_k)\\
&=&\sum_{j,k=1}^m-\frac{\delta_{jk}(1-|\tilde{y}|^2)+\tilde{y}_j\tilde{y}_k}
{(1-|\tilde{y}|^2)^{3/2}}\cdot(x_j-y_j)(x_k-y_k)\\
&=&-\frac{|x-y|^2(1-|\tilde{y}|^2)+\sum_{j,k=1}^m\tilde{y}_j
\tilde{y}_k(x_j-y_j)(x_k-y_k)}{(1-|\tilde{y}|^2)^{3/2}}\\
&=&-\frac{|x-y|^2(1-|\tilde{y}|^2)+|\ll
x-y,\tilde{y}\rl|^2}{(1-|\tilde{y}|^2)^{3/2}},
\end{eqnarray*}
where $\tilde{y}=y+\zeta(x-y)$. Then by \eqref{xy2}, \eqref{xy1} is
proved. \ep

Now we give the proof of Theorem \ref{theoremuabr}.

\noindent{\it Proof of Theorem \ref{theoremuabr}.}\quad  Let $a$,
$b$ and $r$ be fixed. First assume that $b>0$. From Lemma
\ref{lemmaRI}, we have $\l=\l(r,a,b)$ and $\mu=\mu(r,a,b)>0$ such
that $R(r,\l,\mu)=a$ and $I(r,\l,\mu)=b.$ For the need of
\eqref{u4}, let
\begin{equation} \label{u0dy}
u_0(\o)=\frac{A_{r,\l,\mu}(\o)}{\sqrt{1+|A_{r,\l,\mu}(\o)|^2}},
\end{equation}
where $A_{r,\l,\mu}(\o)$ is defined as \eqref{Ady}. Then
$||u_0||_\iy<1$ and by \eqref{RIdy}, we know
\begin{equation} \label{u0td}
\int_S u_0(\o)d\s=R(r,\l,\mu)=a,\quad \int_S \sqrt{1-|u_0(\o)|^2}d\s
=I(r,\l,\mu)=b.
\end{equation}
This means that $u_0\in {\mc U}_{a,b}$.

Let $u\in{\mc U}_{a,b}$. By \eqref{uyq} and \eqref{u0td}, we have
\begin{equation} \label{u1}
\int_S \ll u_0(\o)-u(\o),\l\rl d\s=0,
\end{equation}
\begin{equation} \label{u2}
\mu\int_S (\sqrt{1-|u_0(\o)|^2}-\sqrt{1-|u(\o)|^2})d\s\leq0.
\end{equation}
By Lemma \ref{lemmaxybds}, we have
\begin{equation} \label{u3}
\begin{split}
&\sqrt{1-|u_0(\o)|^2}-\sqrt{1-|u(\o)|^2}\\
&=\frac{\ll u(\o)-u_0(\o),u_0(\o)\rl}{\sqrt{1-|u_0(\o)|^2}}
+\frac{|u(\o)-u_0(\o)|^2(1-|\tilde{u}(\o)|^2)+|\ll
u(\o)-u_0(\o),\tilde{u}(\o)\rl|^2}{2(1-|\tilde{u}(\o)|^2)^{3/2}},
\end{split}\end{equation}
where $\tilde{u}(\o)=u_0(\o)+\zeta(u(\o)-u_0(\o))$, $0<\zeta<1$. By
\eqref{u0dy} and \eqref{Ady}, we have
\begin{equation} \label{u4}
\frac{1}{|rN-\o|^n}l-\l-\frac{\mu u_0(\o)}{\sqrt{1-|u_0(\o)|^2}}=0.
\end{equation}
Then by \eqref{Ldy} and \eqref{u1}-\eqref{u4}, we obtain that
\begin{eqnarray*}
&&\frac{L_r(u_0)-L_r(u)}{1-r^2}=\int_S\frac{\ll
u_0(\o)-u(\o),l\rl}{|rN-\o|^n}d\s\\
&\geq&\int_S\frac{\ll u_0(\o)-u(\o),l\rl}{|rN-\o|^n}d\s-\int_S \ll
u_0(\o)-u(\o),\l\rl
d\s\\
&&\ \ \ \ \ \ \ \ \ \ \ \ \ \ \ \ \
+\mu\int_S (\sqrt{1-|u_0(\o)|^2}-\sqrt{1-|u(\o)|^2})d\s\\
&=&\int_S\ll u_0(\o)-u(\o),\frac{1}{|rN-\o|^n}l-\l-\frac{\mu u_0(\o)}{\sqrt{1-|u_0(\o)|^2}}\rl d\s\\
&&\ \ \ \ \ \ \ \ \ \ \ \ \ \ \ \ \
+\mu\int_S\frac{|u(\o)-u_0(\o)|^2(1-|\tilde{u}(\o)|^2)+|\ll
u(\o)-u_0(\o),\tilde{u}(\o)\rl|^2}{2(1-|\tilde{u}(\o)|^2)^{3/2}}d\s\\
&=&\mu\int_S\frac{|u(\o)-u_0(\o)|^2(1-|\tilde{u}(\o)|^2)+|\ll
u(\o)-u_0(\o),\tilde{u}(\o)\rl|^2}{2(1-|\tilde{u}(\o)|^2)^{3/2}}d\s.
\end{eqnarray*}
Note that
\begin{equation*}
\begin{split}
||\tilde{u}(\o)||&=||u_0(\o)+\zeta(u(\o)-u_0(\o))||\\
&=||u_0(\o)(1-\zeta)+\zeta u(\o)||\\
&\leq||u_0(\o)||(1-\zeta)+||u(\o)||\zeta\\
&<1-\zeta+\zeta=1.
\end{split}
\end{equation*}
Thus $L_r(u_0)\geq L_r(u)$ with equality if and only if
$u(\o)=u_0(\o)$ almost everywhere on $S$. This shows that $u_0(\o)$
is the unique extremal mapping, which will be denoted by
$u_{a,b,r}(\o)$.

Next we consider the case that $b=0$. For the need of \eqref{b=0u5},
let
\begin{equation} \label{b=0u0dy}
u_0(\o)=\frac{\mc A_{r,\l(r,a)}(\o)} {|\mc A_{r,\l(r,a)}(\o)|},
\end{equation}
where $\l(r,a)$ and $\mc A_{r,\l(r,a)}(\o)$ are defined in Lemma
\ref{lemmaR}. Obviously, $||u_0||_\iy\leq1$,
\begin{equation} \label{b=0utd2}
\int_S \sqrt{1-|u_0(\o)|^2}d\s=0,
\end{equation}
and by Lemma \ref{lemmaR},
\begin{equation} \label{b=0u0td}
\int_S u_0(\o)d\s=\mc R(r,\l(r,a))=a.
\end{equation}
This means that $u_0\in {\mc U}_{a,0}$.

Let $u\in{\mc U}_{a,0}$. By \eqref{uyq} and \eqref{b=0u0td}, we have
\begin{equation} \label{b=0u1}
\int_S \ll u_0(\o)-u(\o),\l(r,a)\rl d\s=0.
\end{equation}
By \eqref{A}, we have
\begin{equation} \label{b=0u2}
\frac{1}{|rN-\o|^n}l-\l(r,a)=\mc A_{r,\l(r,a)}(\o).
\end{equation}
By $||u||_\iy\leq1$, we have
\begin{equation} \label{b=0u3}
|\ll u(\o),\mc A_{r,\l(r,a)}(\o)\rl|\leq|u(\o)||\mc
A_{r,\l(r,a)}(\o)|\leq|\mc A_{r,\l(r,a)}(\o)|,
\end{equation}
and
\begin{equation} \label{b=0u4}
|\mc A_{r,\l(r,a)}(\o)|=\ll u(\o),\mc
A_{r,\l(r,a)}(\o)\rl\quad\mbox{if and only if}\quad u(\o)=\frac{\mc
A_{r,\l(r,a)}(\o)} {|\mc A_{r,\l(r,a)}(\o)|}=u_0(\o).
\end{equation}
Then by \eqref{Ldy}, \eqref{b=0u0dy} and
\eqref{b=0u1}-\eqref{b=0u4}, we obtain that
\begin{eqnarray}
\frac{L_r(u_0)-L_r(u)}{1-r^2}
\nonumber &=&\int_S\frac{\ll u_0(\o)-u(\o),l\rl}{|rN-\o|^n}d\s\\
\nonumber &=&\int_S\frac{\ll
u_0(\o)-u(\o),l\rl}{|rN-\o|^n}d\s-\int_S \ll
u_0(\o)-u(\o),\l(r,a)\rl d\s\\
\nonumber &=&\int_S\ll u_0(\o)-u(\o),\frac{1}{|rN-\o|^n}l-\l(r,a)\rl
d\s\\
\nonumber &=&\int_S\ll u_0(\o)-u(\o),\mc A_{r,\l(r,a)}(\o)\rl d\s\\
&=&\int_S\ll \frac{\mc A_{r,\l(r,a)}(\o)}{|\mc
A_{r,\l(r,a)}(\o)|}-u(\o),\mc A_{r,\l(r,a)}(\o)\rl
d\s\label{b=0u5}\\
\nonumber &=&\int_S(|\mc A_{r,\l(r,a)}(\o)|-\ll u(\o),\mc
A_{r,\l(r,a)}(\o)\rl)d\s\geq0
\end{eqnarray}
with equality if and if $u(\o)=u_0(\o)$ almost everywhere on $S$.
Thus $L_r(u_0)\geq L_r(u)$ with equality if and if $u(\o)=u_0(\o)$
almost everywhere on $S$. The theorem is proved. \qed

Let $a\in\mathbb{R}^m$, $b\in\mathbb{R}$, $|a|^2+b^2<1$, and
$0<r<1$. If $b\ge0$, $u_{a,b,r}$ has been defined in Theorem
\ref{theoremuabr}. Now, define
\begin{equation} \label{vdy}
v_{a,b,r}(\o)=\sqrt{1-|u_{a,b,r}(\o)|^2}\quad\mbox{for}\quad\o\in S,
\end{equation}
and
\begin{equation} \label{Udy}
U_{a,b,r}(x)=\int_S\frac{1-|x|^2}{|x-\o|^n} u_{a,b,r}(\o)d\s,
\end{equation}
\begin{equation} \label{Vdy}
V_{a,b,r}(x)=\int_S\frac{1-|x|^2}{|x-\o|^n} v_{a,b,r}(\o)d\s.
\end{equation}
For $b<0$, let
\begin{equation} \label{UVtd}
U_{a,b,r}(x)=U_{a,-b,r}(x),\quad V_{a,b,r}(x)=-V_{a,-b,r}(x).
\end{equation}
Then for any $a\in\mathbb{R}^m$, $b\in\mathbb{R}$ and $|a|^2+b^2<1$,
let
\begin{equation} \label{Fabrdy}
F_{a,b,r}(x)=(U_{a,b,r}(x),V_{a,b,r}(x))\quad\mbox{for}\quad
x\in\mathbb{B}^n.
\end{equation}
The harmonic mapping $F_{a,b,r}(x)$ satisfies $F_{a,b,r}(0)=(a,b)$
and $F_{a,b,r}(\mathbb{B}^n)\subset\mathbb{B}^{m+1}$, since we will
show that $|U_{a,b,r}(x)|^2+|V_{a,b,r}(x)|^2<1$. By the convexity of
the square function,
$$|U_{a,b,r}(x)|^2+|V_{a,b,r}(x)|^2\leq\int_S\frac{1-|x|^2}{|x-\o|^n}
(|u_{a,b,r}(\o)|^2+v_{a,b,r}^2(\o))d\s=1$$ with equality if and only
if $u_{a,b,r,1}(\o),u_{a,b,r,2}(\o),\cdots,u_{a,b,r,m}(\o)$ and
$v_{a,b,r}(\o)$ are constants almost everywhere on $S$, where
$$u_{a,b,r}(\o)=(u_{a,b,r,1}(\o),u_{a,b,r,2}(\o),\cdots,u_{a,b,r,m}(\o)).$$
However $u_{a,b,r,1}(\o),u_{a,b,r,2}(\o),\cdots,u_{a,b,r,m}(\o)$ are
not possiblely constants almost everywhere on $S$. Thus
$|U_{a,b,r}(x)|^2+|V_{a,b,r}(x)|^2<1$.

The mappings $F_{a,b,r}$ are the extremal mappings in the following
theorem. Theorem \ref{theoremUbds} extends Theorem B to
$F\in\O_{n,m+1}$, and when $n=m+1=2$, Theorem \ref{theoremUbds} is
coincident with Theorem B. Note that in the following theorem,
$U_{a,b,r}$ is defined as \eqref{Udy} and \eqref{UVtd}, $F_{a,b,r}$
is defined as \eqref{Fabrdy}.

\bt\label{theoremUbds} Let $F(x)=(U(x),V(x))$ be a harmonic mapping
such that $F(\mathbb{B}^n)\subset\mathbb{B}^{m+1}$ and $F(0)=(a,b)$,
where $U(x)\in\mathbb{R}^{m}$, $V(x)\in\mathbb{R}$,
$a\in\mathbb{R}^{m}$ and $b\in\mathbb{R}$. Let
$l=(1,0,\cdots,0)\in\mathbb{R}^{m}$. Then, for $0<r<1$ and $\o\in
S$,
$$\ll U(r\o),l\rl\le \ll U_{a,b,r}(rN),l\rl$$
with equality at some point $r\o$ if and only if
$F(x)=F_{a,b,r}(xA)$, where $A$ is an orthogonal matrix such that
$\o A=N$. Further, $\ll U(x),l\rl<\ll U_{a,b,r}(rN),l\rl$ for
$|x|<r$. \et

\bp \underline{Step 1}: First the case that $r\o=rN$ will be proved.
Let $0<\tilde{r}<1$ be fixed. Construct mapping
$$G(x)=F(\tilde{r}x)\quad\mbox{for}\quad x\in\ov{\mathbb{B}}^n.$$
$G(x)$ is harmonic on $\ov{\mathbb{B}}^n$ and $G(0)=(a,b)$. Let
$G(x)=(u(x),v(x))$, where $u(x)\in\mathbb{B}^{m}$. Then
$$\|u\|_{\infty}\le 1,\quad\quad \int_S u(\o)d\s=a,$$
\begin{equation} \label{Gtd}
\int_S\sqrt{1-|u(\o)|^2}d\s\ge\int_S |v(\o)|d\s\ge\lt|\int_S
v(\o)d\s\rt|=|b|.
\end{equation}
So by \eqref{uyq} we know that $u\in{\mc U}_{a,|b|}$, and by Theorem
\ref{theoremuabr} we have
$$\ll u(rN),l\rl\le \ll
U_{a,|b|,r}(rN),l\rl$$ with equality if and only if
$u(\o)=u_{a,|b|,r}(\o)$ almost everywhere on $S$. For
$u_{a,|b|,r}(\o)$, by \eqref{u0td} and \eqref{b=0utd2} we have
\begin{equation} \label{jfb}
\int_S \sqrt{1-|u_{a,|b|,r}(\o)|^2}d\s=|b|.
\end{equation}

If $u(\o)=u_{a,|b|,r}(\o)$ almost everywhere on $S$, then by
\eqref{Udy} and \eqref{UVtd}, we have
$$u(x)=U_{a,|b|,r}(x)=U_{a,b,r}(x)\ \ \mbox{for $x\in\mathbb{B}^n$};$$
and by \eqref{vdy}, we have
\begin{equation} \label{vbds}
v_{a,|b|,r}(\o)=\sqrt{1-|u_{a,|b|,r}(\o)|^2}=\sqrt{1-|u(\o)|^2}.
\end{equation}
Note that by \eqref{Gtd}, \eqref{jfb} and \eqref{vbds} we have
$$|b|=\int_S
v_{a,|b|,r}(\o)d\s\geq\int_S|v(\o)|d\s\geq\lt|\int_S
v(\o)d\s\rt|=|b|.$$ Then $$v(\o)=v_{a,|b|,r}(\o)\quad\mbox{almost
everywhere on}\quad S \quad\mbox{when}\quad b\geq0,$$
$$v(\o)=-v_{a,|b|,r}(\o)\quad\mbox{almost
everywhere on}\quad S \quad\mbox{when}\quad b<0.$$ So
$$v(x)=V_{a,b,r}(x)\ \ \mbox{for $x\in\mathbb{B}^n$}.$$

For $G(x)=(u(x),v(x))$, it is proved that $\ll u(rN),l\rl\le \ll
U_{a,b,r}(rN),l\rl$ with equality if and only if
$G(x)=F_{a,b,r}(x)$. Now let $\tilde{r}\rightarrow1$. Note that
$$\lim_{\tilde{r}\rightarrow1}G(x)=\lim_{\tilde{r}\rightarrow1}F(\tilde{r}x)=F(x),
\quad \lim_{\tilde{r}\rightarrow1}u(rN)=U(rN).$$ Then by the result
for $G(x)$, we have $\ll U(rN),l\rl\le \ll U_{a,b,r}(rN),l\rl$ with
equality if and only if $F(x)=F_{a,b,r}(x)$.

\underline{Step 2}: Now we prove the case that $r\o\neq rN$.
Construct mapping
$$\tilde{F}(x)=F(xA^{-1})\quad\mbox{for}\quad x\in\mathbb{B}^n,$$
where $A$ is an orthogonal matrix such that $r\o A=rN$ and $A^{-1}$
is the inverse matrix of $A$. By \cite{ABW}, we know that
$\tilde{F}(x)$ is also a harmonic mapping. Let
$$\tilde{F}(x)=(\tilde{U}(x),\tilde{V}(x)).$$ Note that
$\tilde{F}(0)=F(0)=(a,b)$. Then by the result of step 1, we have
$$\ll\tilde{U}(rN),l\rl\le \ll U_{a,b,r}(rN),l\rl$$ with equality if
and only if $\tilde{F}(x)=F_{a,b,r}(x)$. Note that
$\tilde{U}(rN)=U(rNA^{-1})=U(r\o)$ and $\tilde{F}(x)=F(xA^{-1})$.
Thus $$\ll U(r\o),l\rl\le\ll U_{a,b,r}(rN),l\rl$$ with equality if
and only if $F(xA^{-1})=F_{a,b,r}(x)$. It is just that $\ll
U(r\o),l\rl\le\ll U_{a,b,r}(rN),l\rl$ with equality if and only if
$F(x)=F_{a,b,r}(xA)$.

\underline{Step 3}: We will show that $\ll U(x),l\rl<\ll
U_{a,b,r}(rN),l\rl$ for $|x|<r$. Let
\begin{equation} \label{gdy}
g(x)=\ll U(x),l\rl\ \ \mbox{for $x\in\mathbb{B}^n$}.
\end{equation}
Then $g(x)$ is a real-valued harmonic function. By the result of
step 2, we know that $$g(r\o)\leq\ll U_{a,b,r}(rN),l\rl.$$ Then by
the maximum principle, we have
$$g(x)\leq\ll U_{a,b,r}(rN),l\rl\ \ \mbox{for $|x|\leq r$}.$$

If there exists a point $x_0$ with $|x_0|<r$, such that $g(x_0)=\ll
U_{a,b,r}(rN),l\rl$, then
\begin{equation} \label{gtd}
g(x)\equiv\ll U_{a,b,r}(rN),l\rl\ \ \mbox{for $|x|\leq r$}.
\end{equation}
Then $$g(rN)=\ll U_{a,b,r}(rN),l\rl.$$ Since by \eqref{gdy}
$$g(rN)=\ll U(rN),l\rl,$$ then we have $$\ll U(rN),l\rl=\ll
U_{a,b,r}(rN),l\rl.$$ Then by the result of step 1, we have
$U(x)=U_{a,b,r}(x)$. Thus by \eqref{gdy} and \eqref{gtd}, we obtain
$$\ll U_{a,b,r}(x),l\rl\equiv \ll U_{a,b,r}(rN),l\rl\ \ \mbox{for $|x|\leq
r$}.$$ However, it is impossible since $\ll U_{a,b,r}(x),l\rl$ is
not a constant for $|x|\leq r$. Therefore, for any $x$ with $|x|<r$,
we have $g(x)<\ll U_{a,b,r}(rN),l\rl$. The proof of the theorem is
complete.
 \ep

Consequently, we have a corollary as follows.

\bc\label{e0} Let $F(x)$ be a harmonic mapping such that
$F(\mathbb{B}^n)\subset\mathbb{B}^{m+1}$ and $F(0)=(a,b)$, where
$a\in\mathbb{R}^{m}$ and $b\in\mathbb{R}$. Let
$e_0=(1,0,\cdots,0)\in\mathbb{R}^{m+1}$. Then, for $0<r<1$ and
$\o\in S$,
$$\ll F(r\o),e_0\rl\le \ll F_{a,b,r}(rN),e_0\rl$$
with equality at some point $r\o$ if and only if
$F(x)=F_{a,b,r}(xA)$, where $A$ is an orthogonal matrix such that
$\o A=N$, $F_{a,b,r}$ is defined as \eqref{Fabrdy}. Furthermore,
$\ll F(x),e_0\rl<\ll F_{a,b,r}(rN),e_0\rl$ for $|x|<r$. \ec

Generally, we have Theorem \ref{e} in Section 1. Now we give the
proof of Theorem \ref{e}.

\noindent{\it Proof of Theorem \ref{e}.}\quad  For
$x\in\mathbb{B}^n$, we have $$\ll
F(x),e\rl=F(x)e^T=F(x)(e_0Q_e^{-1})^T=F(x)(e_0Q_e^T)^T
=F(x)Q_ee_0^T=\ll F(x)Q_e,e_0\rl,$$ where $T$ is the transpose
symbol. Let
$$\widetilde{F}(x)=F(x)Q_e,\ \ \ \ x\in\mathbb{B}^n.$$
Then $\widetilde{F}(x)$ is a harmonic mapping by \cite{ABW}, and
$\widetilde{F}(\mathbb{B}^n)\subset\mathbb{B}^{m+1}$,
$\widetilde{F}(0)=F(0)Q_e=(a,b)Q_e$. Using Corollary \ref{e0} to
$\widetilde{F}(x)$, we have for $0<r<1$ and $\o\in S$,
$$\ll \widetilde{F}(r\o),e_0\rl\le \ll F_{(a,b)Q_e,r}(rN),e_0\rl$$
with equality at some point $r\o$ if and only if
$\widetilde{F}(x)=F_{(a,b)Q_e,r}(xA)$, where $A$ is an orthogonal
matrix such that $r\o A=rN$. Furthermore, $\ll
\widetilde{F}(x),e_0\rl<\ll F_{(a,b)Q_e,r}(rN),e_0\rl$ for $|x|<r$.
Note that for $x\in\mathbb{B}^n$, $\widetilde{F}(x)=F(x)Q_e$,
$$\ll
\widetilde{F}(x),e_0\rl=\ll F(x)Q_e,e_0\rl=\ll F(x),e\rl$$ and
$$\ll \widetilde{F}(r\o),e_0\rl=\ll F(r\o),e\rl.$$ Then the theorem
is proved. \qed

From Theorem \ref{e}, we obtain Corollary \ref{especial} in Section
1. Now we give the proof of Corollary \ref{especial}.

\noindent{\it Proof of Corollary \ref{especial}.}\quad  We will
prove the corollary by three steps.

\underline{Step 1}: We claim that for $0<r<1$,
\begin{equation} \label{F0=01}
F_{0,0,r}(x)=(U(x),0,\cdots,0),
\end{equation}
where $U$ is the Poisson integral of the function that equals 1 on
$S^+$ and -1 on $S^-$.

By Theorem \ref{theoremuabr}, \eqref{b=0u0dy}, \eqref{A} and Lemma
\ref{lemmaR}, we have that
\begin{equation*}
u_{0,0,r}(\o)=
\begin{cases}
(1,0,\cdots,0), &\o\in S^+;\\
(-1,0,\cdots,0), &\o\in S^-.
\end{cases}
\end{equation*}
Then by \eqref{Udy}, \eqref{vdy} and \eqref{Vdy}, we obtain that
$$U_{0,0,r}(x)=(U(x),0,\cdots,0)\ \ \mbox{and $V_{0,0,r}(x)\equiv0$}.$$ Thus
$$F_{0,0,r}(x)=(U_{0,0,r}(x),V_{0,0,r}(x))=(U(x),0,\cdots,0).$$ The
claim is proved.

\underline{Step 2}: For any $x\in\mathbb{B}^n$, let $|x|=r, x=r\o$.
Since $F(0)=0$, by Theorem \ref{e}, we have that for
$e_0=(1,0,\cdots,0)\in\mathbb{R}^{m+1}$ and any unit vector
$e\in\mathbb{R}^{m+1}$,
$$\ll F(r\o),e\rl\le \ll F_{0,0,r}(rN),e_0\rl.$$
That is
\begin{equation} \label{F0=02}
\ll F(x),e\rl\le \ll F_{0,0,|x|}(|x|N),e_0\rl.
\end{equation}
If $F(x)=0$, then obviously $|F(x)|\leq U(|x|N)$ since
$U(|x|N)\geq0$. If $F(x)\neq0$, then let $e=\frac{F(x)}{|F(x)|}$ and
consequently by \eqref{F0=01} and \eqref{F0=02}, we have $|F(x)|\leq
U(|x|N)$.

\underline{Step 3}: For some $x_0\in\mathbb{B}^n$, let $|x_0|=r_0$.
By Step 2 and Theorem \ref{e}, we have that $|F(x_0)|=U(|x_0|N)$ if
and only if $F(x)=F_{0,0,r_0}(xA)Q_e^{-1}$, where $A$ is an
orthogonal matrix such that $x_0A=r_0N$,
$e=\frac{F(x_0)}{|F(x_0)|}$, $Q_e$ be an orthogonal matrix such that
$eQ_e=e_0$, $Q_e^{-1}$ is the inverse matrix of $Q_e$. By
\eqref{F0=01},
$$F_{0,0,r_0}(xA)=(U(xA),0,\cdots,0).$$ Note that
$$(U(xA),0,\cdots,0)=(U(xA),0,\cdots,0)e_0^Te_0,$$ where $T$ is the
transpose symbol. Then $$F(x)=(U(xA),0,\cdots,0)Q_e^{-1}
=((U(xA),0,\cdots,0)e_0^T)(e_0Q_e^{-1})=U(xA)e.$$ The corollary is
proved.\qed

\section{The proofs of Lemma 1 and Lemma 2}\label{The proofs}
For the proofs of Lemma 1 and Lemma 2, we need the following two
lemmas.

\bl\label{lemmaQn} Let the determinant
\begin{equation*}
Q_n=\lt|
\begin{array}{cccc}
b+a_{11}&a_{12}&\cdots&a_{1n}\\
a_{21}&b+a_{22}&\cdots&a_{2n}\\
\vdots&\vdots&\ddots&\vdots\\
a_{n1}&a_{n2}&\cdots&b+a_{nn},\end{array} \rt| \ \ \ \ \ n\geq2,
\end{equation*}
where $a_{ij}=-c_ic_j$ for $i\neq1$ or $j\neq1$. Then
\begin{equation}\label{Qn}
Q_n=b^n+b^{n-1}(a_{11}+a_{22}+\cdots+a_{nn})+b^{n-2}\lt( \lt|
\begin{array}{cc}
a_{11}&a_{12}\\
a_{21}&a_{22}\end{array}\rt|+\lt|
\begin{array}{cc}
a_{11}&a_{13}\\
a_{31}&a_{33}\end{array}\rt|+\cdots+\lt|
\begin{array}{cc}
a_{11}&a_{1n}\\
a_{n1}&a_{nn}\end{array}\rt|\rt).
\end{equation}\el

\bp Let \begin{equation*} A_n=\lt|
\begin{array}{cccc}
a_{11}&a_{12}&\cdots&a_{1n}\\
a_{21}&a_{22}&\cdots&a_{2n}\\
\multicolumn{4}{c}\dotfill\\
a_{n1}&a_{n2}&\cdots&a_{nn}\end{array}\rt|.
\end{equation*}
By the definition of the determinant $Q_n$, we have
\begin{equation}\label{Qn1}
\begin{split}
Q_n=b^n&+b^{n-1}(\mbox{the sum of all the level 1 principal minor of
$A_n$})\\
&+b^{n-2}(\mbox{the sum of all the level 2 principal minor of $
A_n$)}+\cdots\\
&+b(\mbox{the sum of all the level n-1 principal minor of
$A_n$)}+A_n.
\end{split}
\end{equation}

When $n=2$, \eqref{Qn} obviously holds.

When $n\geq3$, we will prove that for integer $k\geq3$, the value of
any level $k$ principal minor of $A_n$ is $0$. Let $P_k$ is a level
$k$ principal minor of $A_n$ with $k\geq3$. Denote
\begin{equation*} P_k=\lt|
\begin{array}{cccc}
a_{i_1i_1}&a_{i_1i_2}&\cdots&a_{i_1i_k}\\
a_{i_2i_1}&a_{i_2i_2}&\cdots&a_{i_2i_k}\\
\multicolumn{4}{c}\dotfill\\
a_{i_ki_1}&a_{i_ki_2}&\cdots&a_{i_ki_k}\end{array}\rt|.
\end{equation*}
Let
\begin{equation*}
\begin{split}
A((k-1)k|pq)=&\mbox{the $2\times2$ minor
 of $P_k$ that lies on the intersection of rows $(k-1),k$}\\
 &\mbox{with columns $p,q$,
  where $1\leq p<q\leq k$},
\end{split}
\end{equation*}
and
\begin{equation*}
\begin{split}
M((k-1)k|pq)=&\mbox{the $(n-2)\times (n-2)$ minor obtained by
deleting rows $(k-1),k$}\\
 &\mbox{and columns $p,q$ from $P_k$,
  where $1\leq p<q\leq k$}.
\end{split}
\end{equation*}
The cofactor of $A((k-1)k|pq)$ is defined to be the signed minor
$$\tilde{A}((k-1)k|pq)=(-i)^{((k-1)+k+p+q)}M((k-1)k|pq).$$
Then using Laplace's expansion to evaluate $P_k$ in terms of the
last two rows, we have
\begin{equation}\label{Pk}
P_k=\sum_{1\leq p<q\leq k}A((k-1)k|pq)\tilde{A}((k-1)k|pq).
\end{equation}
By $k-1\geq2$, we have for any $p,q\ (1\leq p<q\leq k)$,
\begin{equation}\label{Pkz}
\begin{split}
A((k-1)k|pq)&=\lt|
\begin{array}{cc}
a_{i_{k-1}i_p}&a_{i_{k-1}i_q}\\
a_{i_ki_p}&a_{i_ki_q}\end{array}\rt|\\
&=\lt|
\begin{array}{cc}
-c_{i_{k-1}}c_{i_p}&-c_{i_{k-1}}c_{i_q}\\
-c_{i_{k}}c_{i_p}&-c_{i_{k}}c_{i_q}\end{array}\rt|\\
&=0.
\end{split}
\end{equation}
Thus $P_k=0$ by \eqref{Pk} and \eqref{Pkz}.

When $n\geq3$, it is proved above that for integer $k\geq3$, the
value of any level $k$ principal minor of $A_n$ is $0$. Then by
\eqref{Qn1},
\begin{equation*}
\begin{split}
Q_n&=b^n+b^{n-1}(\mbox{the sum of all the level 1 principal minor of
$A_n$})\\
&\ \ \ \ \ \ \ \ +b^{n-2}(\mbox{the sum of all the level 2 principal
minor
of $ A_n$)}\\
&=b^n+b^{n-1}(a_{11}+a_{22}+\cdots+a_{nn})+b^{n-2}\sum_{1\leq
i<j\leq n}\lt|
\begin{array}{cc}
a_{ii}&a_{ij}\\
a_{ji}&a_{jj}\end{array}\rt|.
\end{split}
\end{equation*}
Note that
\begin{equation*}
\lt|
\begin{array}{cc}
a_{ii}&a_{ij}\\
a_{ji}&a_{jj}\end{array}\rt|=\lt|
\begin{array}{cc}
-c_ic_i&-c_ic_j\\
-c_jc_i&-c_jc_j\end{array}\rt|=0\ \ \ \ \ \mbox{for $1<i<j\leq n$}.
\end{equation*}
Thus \eqref{Qn} holds when $n\geq3$. Then the lemma is proved.
 \ep

\bl\label{lemmafcz} Fixed integer $k\geq1$, let matrices
\begin{equation*} A=\lt(a_{ij}\rt)_{k\times k},\ \ x=\lt(
\begin{array}{c}
x_1\\
x_2\\
\vdots\\
x_k\end{array}\rt),\ \ b=\lt(
\begin{array}{c}
b_1\\
b_2\\
\vdots\\
b_k\end{array}\rt),\ \ c=(c_1,c_2,\cdots,c_k),\ \ B=\lt(
\begin{array}{cc}
A&b\\
c&c_{k+1}
\end{array}\rt).
\end{equation*}
Suppose that $Ax+b=0$ and $det(A)\neq0$. Then
\begin{equation}\label{c}
cx+c_{k+1}=\frac{det(B)}{det(A)}.
\end{equation}
\el

\bp For $1\leq i\leq k$, let $A_i$ be the determinant obtained by
replacing column $j$ of $det(A)$ with $-b$. For $1\leq j\leq k+1$,
let $B_j$ be the determinant obtained by deleting row $k+1$ and
column $j$ from $det(B)$.

Using Cramer's rule to $Ax+b=0$ we have
$$x_1=\frac{A_1}{det(A)},x_2=\frac{A_2}{det(A)},\cdots,x_k=\frac{A_k}{det(A)}.$$
Then
\begin{equation}\label{c1}
c_1x_1+c_2x_2+\cdots+c_kx_k+c_{k+1}=\frac{1}{det(A)}(c_1A_1+c_2A_2+\cdots+c_kA_k+c_{k+1}det(A)).
\end{equation}
Note that
\begin{equation}\label{c2}
\begin{split}
&\ \ c_1A_1+c_2A_2+\cdots+c_kA_k+c_{k+1}det(A)\\
&=c_1\lt|
\begin{array}{ccccc}
-b_1&a_{12}&a_{13}&\cdots&a_{1k}\\
-b_2&a_{22}&a_{23}&\cdots&a_{2k}\\
\multicolumn{5}{c}\dotfill\\
-b_k&a_{k2}&a_{k3}&\cdots&a_{kk}\end{array}\rt|+c_2\lt|
\begin{array}{ccccc}
a_{11}&-b_1&a_{13}&\cdots&a_{1k}\\
a_{21}&-b_2&a_{23}&\cdots&a_{2k}\\
\multicolumn{5}{c}\dotfill\\
a_{k1}&-b_k&a_{k3}&\cdots&a_{kk}\end{array}\rt|+\cdots\\
&\ \ \ \ \ \ \ \ +c_k\lt|
\begin{array}{ccccc}
a_{11}&a_{12}&\cdots&a_{1(k-1)}&-b_1\\
a_{21}&a_{22}&\cdots&a_{2(k-1)}&-b_2\\
\multicolumn{5}{c}\dotfill\\
a_{k1}&a_{k2}&\cdots&a_{k(k-1)}&-b_k\end{array}\rt|+c_{k+1}\lt|
\begin{array}{cccc}
a_{11}&a_{12}&\cdots&a_{1k}\\
a_{21}&a_{22}&\cdots&a_{2k}\\
\multicolumn{4}{c}\dotfill\\
a_{k1}&a_{k2}&\cdots&a_{kk}\end{array}\rt|\\
&=c_1(-1)^{k}\lt|
\begin{array}{ccccc}
a_{12}&a_{13}&\cdots&a_{1k}&b_1\\
a_{22}&a_{23}&\cdots&a_{2k}&b_2\\
\multicolumn{5}{c}\dotfill\\
a_{k2}&a_{k3}&\cdots&a_{kk}&b_k\end{array}\rt|+c_2(-1)^{k-1}\lt|
\begin{array}{ccccc}
a_{11}&a_{13}&\cdots&a_{1k}&b_1\\
a_{21}&a_{23}&\cdots&a_{2k}&b_2\\
\multicolumn{5}{c}\dotfill\\
a_{k1}&a_{k3}&\cdots&a_{kk}&b_k\end{array}\rt|+\cdots\\
&\ \ \ \ \ \ \ \ +c_k(-1)^{1}\lt|
\begin{array}{ccccc}
a_{11}&a_{12}&\cdots&a_{1(k-1)}&b_1\\
a_{21}&a_{22}&\cdots&a_{2(k-1)}&b_2\\
\multicolumn{5}{c}\dotfill\\
a_{k1}&a_{k2}&\cdots&a_{k(k-1)}&b_k\end{array}\rt|+c_{k+1}(-1)^{0}\lt|
\begin{array}{cccc}
a_{11}&a_{12}&\cdots&a_{1k}\\
a_{21}&a_{22}&\cdots&a_{2k}\\
\multicolumn{4}{c}\dotfill\\
a_{k1}&a_{k2}&\cdots&a_{kk}\end{array}\rt|\\
&=c_1(-1)^{k}B_1+c_2(-1)^{k-1}B_2+\cdots+c_k(-1)^{1}B_k+c_{k+1}(-1)^{0}B_{k+1}\\
&=c_1(-1)^{k+2}B_1+c_2(-1)^{k+3}B_2+\cdots+c_k(-1)^{2k+1}B_k+c_{k+1}(-1)^{2k+2}B_{k+1}\\
&=det(B).
\end{split}
\end{equation}
Thus by \eqref{c1} and \eqref{c2}, \eqref{c} is proved. \ep

Now we give the proof of Lemma 1.

\noindent{\it Proof of Lemma 1.}\quad  We will prove Lemma 1 by six
steps, where Step 2 is only for the case that $m=1$, and Step 3 -
Step 5 are only for the case that $m\geq2$.

\underline{Step 1}: We give some denotation and calculation. Write
$$A_{r,\l,\mu}(\omega)=A(\o)=(A_1(\o),A_2(\o),\cdots,A_m(\o)),$$
$$R(r,\l,\mu)=(R_1(r,\l,\mu),R_2(r,\l,\mu),\cdots,R_m(r,\l,\mu)),$$
$$l=(l_1,\cdots,l_m), \l=(\l_1,\l_2,\cdots,\l_m), \ \mbox{and} \ \
a=(a_1,a_2,\cdots,a_m).$$ For $i,j=1,2,\cdots,m$, we denote
$\frac{\p R_j(r,\l,\mu)}{\p\l_i}=R_{ji}$, $\frac{\p
R_j(r,\l,\mu)}{\p\mu}=R_{j\mu}$, $\frac{\p I(r,\l,\mu)}{\p\l_j}=I_j$
and $\frac{\p I(r,\l,\mu)}{\p\mu}=I_\mu$. Then a simple calculation
gives
\begin{equation}\label{Rjjdy}
R_{jj}=-\frac1{\mu}\int_S
\frac{1+|A(\o)|^2-A^2_j(\o)}{(1+|A(\o)|^2)^{3/2}}d\s\ \ \ \
\mbox{for $j=1,2,\cdots,m$};
\end{equation}
\begin{equation}\label{Rjidy}
R_{ji}=-\frac1{\mu}\int_S
\frac{-A_i(\o)A_j(\o)}{(1+|A(\o)|^2)^{3/2}}d\s\ \ \ \ \mbox{for
$i\neq j; i,j=1,2,\cdots,m$};
\end{equation}
\begin{equation}
R_{j\mu}=-\frac1{\mu}\int_S \frac{A_j(\o)}{(1+|A(\o)|^2)^{3/2}}d\s\
\ \ \ \mbox{for $j=1,2,\cdots,m$};
\end{equation}
\begin{equation}
I_j=\frac1{\mu}\int_S \frac{A_j(\o)}{(1+|A(\o)|^2)^{3/2}}d\s\ \ \ \
\mbox{for $j=1,2,\cdots,m$};
\end{equation}
\begin{equation}\label{Imudy}
I_\mu=\frac1{\mu}\int_S \frac{|A(\o)|^2}{(1+|A(\o)|^2)^{3/2}}d\s.
\end{equation}
It is easy to see that\\
(i) by \eqref{Rjjdy}, for $j=1,2,\cdots,m$, $R_{jj}<0$ for any
$\l\in\mathbb{R}^m$ and $\mu>0$, and $R_j(r,\l,\mu)$ is strictly
decreasing as a function of $\l_j$ for fixed the other components of
$\l$ and $\mu$;\\
(ii) by \eqref{Ady} and \eqref{RIdy}, for $j=1,2,\cdots,m$, fixing
$\mu$ and the components of $\l$ expect $\l_j$, $R_j(r,\l,\mu)\to-1$
or $1$ according to $\l_j\to+\iy$ or $\l_j\to-\iy$;\\
(iii) by \eqref{Ady} and \eqref{RIdy}, $0<I(r,\l,\mu)<1$ for any
$\l\in\mathbb{R}^m$ and $\mu>0$.

In addition, we claim that
\begin{equation}\label{1}
\frac{\lt|
\begin{array}{cccc}
R_{11}&R_{12}&\cdots&R_{1(k+1)}\\
R_{21}&R_{22}&\cdots&R_{2(k+1)}\\
\multicolumn{4}{c}\dotfill\\
R_{(k+1)1}&R_{(k+1)2}&\cdots&R_{(k+1)(k+1)}\end{array}\rt|}{\lt|
\begin{array}{cccc}
R_{11}&R_{12}&\cdots&R_{1k}\\
R_{21}&R_{22}&\cdots&R_{2k}\\
\multicolumn{4}{c}\dotfill\\
R_{k1}&R_{k2}&\cdots&R_{kk}\end{array}\rt|}<0\ \ \ \mbox{for integer
$k$ with $1\leq k\leq m-1$ when $m\geq2$};
\end{equation}
and
\begin{equation}\label{2}
\frac{\lt|
\begin{array}{ccccc}
R_{11}&R_{12}&\cdots&R_{1m}&R_{1\mu}\\
R_{21}&R_{22}&\cdots&R_{2m}&R_{2\mu}\\
\multicolumn{4}{c}\dotfill\\
R_{m1}&R_{m2}&\cdots&R_{mm}&R_{m\mu}\\
I_1&I_2&\cdots&I_m&I_\mu\end{array}\rt|}{\lt|
\begin{array}{cccc}
R_{11}&R_{12}&\cdots&R_{1m}\\
R_{21}&R_{22}&\cdots&R_{2m}\\
\multicolumn{4}{c}\dotfill\\
R_{m1}&R_{m2}&\cdots&R_{mm}\end{array}\rt|}>0\ \ \ \mbox{when
$m\geq1$}.
\end{equation}
Now we will prove the two claims above.

For \eqref{Rjjdy}-\eqref{Imudy}, let
$d\tilde{\s}=(1/(1+|A(\o)|^2)^{3/2})d\s$, $T=\int_S d\tilde{\s}$,
$d\xi=(1/T)d\tilde{\s}$, $\tilde{b}=\int_S (1+|A(\o)|^2)d\xi$, and
for $i,j=1,2,\cdots,m$, $\tilde{a}_{ij}=\int_S -A_i(\o)A_j(\o)d\xi$,
$c_j=\int_S A_j(\o)d\xi$. Then $T>0$, $\int_S d\xi=1$, and
\begin{equation}\label{Rjjdy2}
R_{jj}=-\frac{T}{\mu}(\tilde{b}+\tilde{a}_{jj})\ \ \ \ \mbox{for
$j=1,2,\cdots,m$};
\end{equation}
\begin{equation}\label{Rjidy2}
R_{ji}=-\frac{T}{\mu}\tilde{a}_{ij}\ \ \ \ \mbox{for $i\neq j;
i,j=1,2,\cdots,m$};
\end{equation}
\begin{equation}\label{Rjmudy2}
R_{j\mu}=-\frac{T}{\mu}c_j\ \ \ \ \mbox{for $j=1,2,\cdots,m$};
\end{equation}
\begin{equation}\label{Ijdy2}
I_j=\frac{T}{\mu}c_j\ \ \ \ \mbox{for $j=1,2,\cdots,m$};
\end{equation}
\begin{equation}\label{Imudy2}
I_\mu=\frac{T}{\mu}(\tilde{b}-1);
\end{equation}
\begin{equation}\label{bja}
\tilde{b}+\tilde{a}_{11}+\tilde{a}_{22}+\cdots+\tilde{a}_{jj}
\geq\tilde{b}+\tilde{a}_{11}+\tilde{a}_{22}+\cdots+\tilde{a}_{mm}=1\
\ \mbox{for $j=1,2,\cdots,m$}.
\end{equation}
Since $A_1(\o)=\frac{1}{\mu}\lt(\frac{1}{|rN-\o|^n}-\l_1\rt)$ by
\eqref{Ady} and $\int_S d\xi=1$, we have
\begin{equation}\label{a11c1}
\begin{split}
-\tilde{a}_{11}-c_1^2&=\int_S A_1^2(\o)d\xi-\lt(\int_S
A_1(\o)d\xi\rt)^2\\
&=\int_S \lt[A_1(\o)-\int_S A_1(\o)d\xi\rt]^2d\xi>0.
\end{split}
\end{equation}
When $m\geq2$, since $\int_S d\xi=1$ and $A_j(\o)=-\frac{\l_j}{\mu}$
for $j=2,\cdots,m$ by \eqref{Ady}, we have
\begin{equation}\label{aij}
\tilde{a}_{ij}=-\int_S A_i(\o)d\xi\int_S A_j(\o)d\xi=-c_ic_j\ \ \ \
\mbox{for $i\neq1$ or $j\neq1$}, i,j=1,2,\cdots,m;
\end{equation}
and by \eqref{a11c1}, we have
\begin{equation}\label{xhls}
\lt|\begin{array}{cc}
\tilde{a}_{11}&\tilde{a}_{1j}\\
\tilde{a}_{j1}&\tilde{a}_{jj}
\end{array}\rt|=\lt|\begin{array}{cc}
\tilde{a}_{11}&-c_1c_j\\
-c_jc_1&-c_jc_j
\end{array}\rt|
=c_j^2(-\tilde{a}_{11}-c_1^2)\geq0\ \ \ \mbox{for $j=2,\cdots,m$}.
\end{equation}

For integer $1\leq p\leq m$, let
\begin{equation}\label{Qp}
Q_p=\lt|
\begin{array}{cccc}
\tilde{b}+\tilde{a}_{11}&\tilde{a}_{12}&\cdots&\tilde{a}_{1p}\\
\tilde{a}_{21}&\tilde{b}+\tilde{a}_{22}&\cdots&\tilde{a}_{2p}\\
\vdots&\vdots&\ddots&\vdots\\
\tilde{a}_{p1}&\tilde{a}_{p2}&\cdots&\tilde{b}+\tilde{a}_{pp}\end{array}
\rt|.
\end{equation}
By \eqref{bja}, we have that when $p=1$,
$Q_1=\tilde{b}+\tilde{a}_{11}>0$. By \eqref{aij} and Lemma
\ref{lemmaQn}, we have that when $p\geq2$,
\begin{equation*}
\begin{split}
Q_p&=\tilde{b}^p+\tilde{b}^{p-1}\sum_{j=1}^p\tilde{a}_{jj}
+\tilde{b}^{p-2}\sum_{j=2}^p\lt|
\begin{array}{cc}
\tilde{a}_{11}&\tilde{a}_{1j}\\
\tilde{a}_{j1}&\tilde{a}_{jj}\end{array}\rt|\\
&=\tilde{b}^{p-1}\lt(\tilde{b}+\sum_{j=1}^p\tilde{a}_{jj}\rt)
+\tilde{b}^{p-2}\sum_{j=2}^p \lt|
\begin{array}{cc}
\tilde{a}_{11}&\tilde{a}_{1j}\\
\tilde{a}_{j1}&\tilde{a}_{jj}\end{array}\rt|.
\end{split}
\end{equation*}
Consequently, by \eqref{bja}, \eqref{xhls} and $\tilde{b}>0$, we
obtain that when $p\geq2$, $Q_p>0$.

Let
\begin{equation}\label{Qm+1}
Q_{m+1}=\lt|
\begin{array}{ccccc}
\tilde{b}+\tilde{a}_{11}&\tilde{a}_{12}&\cdots&\tilde{a}_{1m}&-c_1\\
\tilde{a}_{21}&\tilde{b}+\tilde{a}_{22}&\cdots&\tilde{a}_{2m}&-c_2\\
\vdots&\vdots&\ddots&\vdots&\vdots\\
\tilde{a}_{m1}&\tilde{a}_{m2}&\cdots&\tilde{b}+\tilde{a}_{mm}&-c_m\\
-c_1&-c_2&\cdots&-c_m&\tilde{b}-1
\end{array}
\rt|.
\end{equation}
If $m=1$, then
\begin{equation*}
Q_{m+1}=\lt|
\begin{array}{cc}
\tilde{b}+\tilde{a}_{11}&-c_1\\
-c_1&\tilde{b}-1
\end{array}
\rt|=\tilde{b}(\tilde{b}+\tilde{a}_{11}-1)+(-\tilde{a}_{11}-c_1^2),
\end{equation*}
and $Q_{m+1}>0$ since \eqref{bja} and \eqref{a11c1}. If $m\geq2$,
then by \eqref{aij}, $-c_j=-c_j\times1, -1=-1\times1$, and Lemma
\ref{lemmaQn}, we have
\begin{equation*}
\begin{split}
Q_{m+1}&=\tilde{b}^{m+1}+\tilde{b}^{m}\lt(\sum_{j=1}^m\tilde{a}_{jj}-1\rt)
+\tilde{b}^{m-1}\lt(\sum_{j=2}^{m} \lt|
\begin{array}{cc}
\tilde{a}_{11}&\tilde{a}_{1j}\\
\tilde{a}_{j1}&\tilde{a}_{jj}\end{array}\rt|+\lt|
\begin{array}{cc}
\tilde{a}_{11}&-c_1\\
-c_1&-1\end{array}\rt|\rt)\\
&=\tilde{b}^{m}\lt(\tilde{b}+\sum_{j=1}^m\tilde{a}_{jj}-1\rt)
+\tilde{b}^{m-1}\lt(\sum_{j=2}^{m} \lt|
\begin{array}{cc}
\tilde{a}_{11}&\tilde{a}_{1j}\\
\tilde{a}_{j1}&\tilde{a}_{jj}\end{array}\rt|+\lt|
\begin{array}{cc}
\tilde{a}_{11}&-c_1\\
-c_1&-1\end{array}\rt|\rt),
\end{split}
\end{equation*}
and  $Q_{m+1}>0$ since \eqref{bja}, \eqref{a11c1}, \eqref{xhls} and
$\tilde{b}>0$.

By \eqref{Rjjdy2}, \eqref{Rjidy2} and \eqref{Qp}, we have for
integer $k$ with $1\leq k\leq m-1$ when $m\geq2$,
\begin{equation*}
\begin{split}
\frac{\lt|
\begin{array}{cccc}
R_{11}&R_{12}&\cdots&R_{1(k+1)}\\
R_{21}&R_{22}&\cdots&R_{2(k+1)}\\
\multicolumn{4}{c}\dotfill\\
R_{(k+1)1}&R_{(k+1)2}&\cdots&R_{(k+1)(k+1)}\end{array}\rt|}{\lt|
\begin{array}{cccc}
R_{11}&R_{12}&\cdots&R_{1k}\\
R_{21}&R_{22}&\cdots&R_{2k}\\
\multicolumn{4}{c}\dotfill\\
R_{k1}&R_{k2}&\cdots&R_{kk}\end{array}\rt|}
&=\frac{\lt(-\frac{T}{\mu}\rt)^{k+1} \lt|
\begin{array}{cccc}
\tilde{b}+\tilde{a}_{11}&\tilde{a}_{12}&\cdots&\tilde{a}_{1(k+1)}\\
\tilde{a}_{21}&\tilde{b}+\tilde{a}_{22}&\cdots&\tilde{a}_{2(k+1)}\\
\vdots&\vdots&\ddots&\vdots\\
\tilde{a}_{(k+1)1}&\tilde{a}_{(k+1)2}&\cdots&\tilde{b}+\tilde{a}_{(k+1)(k+1)}\end{array}
\rt|}{\lt(-\frac{T}{\mu}\rt)^k\lt|
\begin{array}{cccc}
\tilde{b}+\tilde{a}_{11}&\tilde{a}_{12}&\cdots&\tilde{a}_{1k}\\
\tilde{a}_{21}&\tilde{b}+\tilde{a}_{22}&\cdots&\tilde{a}_{2k}\\
\vdots&\vdots&\ddots&\vdots\\
\tilde{a}_{k1}&\tilde{a}_{k2}&\cdots&\tilde{b}+\tilde{a}_{kk}\end{array}
\rt|}\\
&=
\frac{\lt(-\frac{T}{\mu}\rt)^{k+1}Q_{k+1}}{\lt(-\frac{T}{\mu}\rt)^{k}Q_k}\\
&=\lt(-\frac{T}{\mu}\rt)\frac{Q_{k+1}}{Q_k}.
\end{split}\end{equation*}
Note that $T>0, \mu>0, Q_k>0, Q_{k+1}>0$. Then the first claim
\eqref{1} is proved.

By \eqref{Rjjdy2}-\eqref{Imudy2}, \eqref{Qp} and \eqref{Qm+1}, we
have when $m\geq1$,
\begin{equation*}
\begin{split}
\frac{\lt|
\begin{array}{ccccc}
R_{11}&R_{12}&\cdots&R_{1m}&R_{1\mu}\\
R_{21}&R_{22}&\cdots&R_{2m}&R_{2\mu}\\
\multicolumn{4}{c}\dotfill\\
R_{m1}&R_{m2}&\cdots&R_{mm}&R_{m\mu}\\
I_1&I_2&\cdots&I_m&I_\mu\end{array}\rt|}{\lt|
\begin{array}{cccc}
R_{11}&R_{12}&\cdots&R_{1m}\\
R_{21}&R_{22}&\cdots&R_{2m}\\
\multicolumn{4}{c}\dotfill\\
R_{m1}&R_{m2}&\cdots&R_{mm}\end{array}\rt|}
&=\frac{(-1)^m\lt(\frac{T}{\mu}\rt)^{m+1} \lt|
\begin{array}{ccccc}
\tilde{b}+\tilde{a}_{11}&\tilde{a}_{12}&\cdots&\tilde{a}_{1m}&-c_1\\
\tilde{a}_{21}&\tilde{b}+\tilde{a}_{22}&\cdots&\tilde{a}_{2m}&-c_2\\
\vdots&\vdots&\ddots&\vdots&\vdots\\
\tilde{a}_{m1}&\tilde{a}_{m2}&\cdots&\tilde{b}+\tilde{a}_{mm}&-c_m\\
-c_1&-c_2&\cdots&-c_m&\tilde{b}-1
\end{array} \rt|}{\lt(-\frac{T}{\mu}\rt)^m \lt|
\begin{array}{cccc}
\tilde{b}+\tilde{a}_{11}&\tilde{a}_{12}&\cdots&\tilde{a}_{1m}\\
\tilde{a}_{21}&\tilde{b}+\tilde{a}_{22}&\cdots&\tilde{a}_{2m}\\
\vdots&\vdots&\ddots&\vdots\\
\tilde{a}_{m1}&\tilde{a}_{m2}&\cdots&\tilde{b}+\tilde{a}_{mm}\end{array}
\rt|}\\
&=
\frac{(-1)^m\lt(\frac{T}{\mu}\rt)^{m+1}Q_{m+1}}{\lt(-\frac{T}{\mu}\rt)^{m}Q_m}\\
&=\frac{T}{\mu}\frac{Q_{m+1}}{Q_m}.
\end{split}\end{equation*}
Note that $T>0, \mu>0, Q_m>0, Q_{m+1}>0$. Then the second claim
\eqref{2} is proved.

\underline{Step 2}: Step 2 is only for the case that $m=1$. By (i)
and (ii) in Step 1, we know that for fixed $\mu$, $R_(r,\l,\mu)$ is
strictly decreasing from $1$ to $-1$ as $\l$ increasing from $-\iy$
to $+\iy$. Then for any $-1<a<1$ and fixed $\mu$, there exists a
unique real number $\l(\mu,a)$ such that
\begin{equation*}
R(r,\l,\mu)\lt|_{\l=\l(\mu,a)}\rt.=a.
\end{equation*}
Further, using the implicit function theorem, we have that the
function $\l=\l(\mu,a)$ defined on $\{(\mu,a):\mu>0,-1<a<1\}$ is a
continuous function and $ \frac{\p\l(\mu,a)}{\p\mu}$ exist.

\underline{Step 3}: Step 3 is only for the case that $m\geq2$. By
(i) and (ii) in Step 1, we know that for fixed $\l_2,\cdots,\l_m$
and $\mu$, $R_1(r,\l,\mu)$ is strictly decreasing from $1$ to $-1$
as $\l_1$ increasing from $-\iy$ to $+\iy$. Then for any $-1<a_1<1$
and fixed $\l_2,\cdots,\l_m$ and $\mu$, there exists a unique real
number $\l_1(\l_2,\cdots,\l_m,\mu,a_1)$ such that
\begin{equation*}
R_1(r,\l,\mu)\lt|_{\l_1=\l_1(\l_2,\cdots,\l_m,\mu,a_1)}\rt.=a_1.
\end{equation*}
Further, using the implicit function theorem, we have that the
function $$\l_1=\l_1(\l_2,\cdots,\l_m,\mu,a_1)$$ defined on
$\{(\l_2,\cdots,\l_m,\mu,a_1):\l_2\in\mathbb{R},\cdots,\l_m\in\mathbb{R},\mu>0,-1<a_1<1\}$
is a continuous function and
$\frac{\p\l_1(\l_2,\cdots,\l_m,\mu,a_1)}{\p\l_2},
\cdots,\frac{\p\l_1(\l_2,\cdots,\l_m,\mu,a_1)}{\p\l_m},
\frac{\p\l_1(\l_2,\cdots,\l_m,\mu,a_1)}{\p\mu}$ exist.

\underline{Step 4}: For the case that $m\geq2$, we will prove the
following result:

For an integer $k$ with $1\leq k\leq m-1$, if\\
(1) there exists a unique continuous function
$\l_1=\l_1(\l_2,\cdots,\l_m,\mu,a_1)$, which defined on\\
$\{(\l_2,\cdots,\l_m,\mu,a_1):\l_2\in\mathbb{R},\cdots,\l_m\in\mathbb{R},\mu>0,-1<a_1<1\}$,
such that
\begin{equation*}
R_1(r,\l,\mu)\lt|_{\l_1=\l_1(\l_2,\cdots,\l_m,\mu,a_1)}\rt.=a_1,
\end{equation*}
and $\frac{\p\l_1(\l_2,\cdots,\l_m,\mu,a_1)}{\p\l_2},
\cdots,\frac{\p\l_1(\l_2,\cdots,\l_m,\mu,a_1)}{\p\l_m},
\frac{\p\l_1(\l_2,\cdots,\l_m,\mu,a_1)}{\p\mu}$ exist;\\
(2) there exists a unique continuous function
$\l_2=\l_2(\l_3,\cdots,\l_m,\mu,a_1,a_2)$, which defined on
$\{(\l_3,\cdots,\l_m,\mu,a_1,a_2):\l_3\in\mathbb{R},\cdots,\l_m\in\mathbb{R},
\mu>0,a_1\in\mathbb{R},a_2\in\mathbb{R},a_1^2+a_2^2<1\}$, such that
\begin{equation*}
R_2(r,\l,\mu)\lt|_{\l_1=\l_1(\l_2,\cdots,\l_m,\mu,a_1)
\atop\l_2=\l_2(\l_3,\cdots,\l_m,\mu,a_1,a_2)}\rt.=a_2,
\end{equation*}
and $\frac{\p\l_2(\l_3,\cdots,\l_m,\mu,a_1,a_2)}{\p\l_3},
\cdots,\frac{\p\l_2(\l_3,\cdots,\l_m,\mu,a_1,a_2)}{\p\l_m},
\frac{\p\l_2(\l_3,\cdots,\l_m,\mu,a_1,a_2)}{\p\mu}$ exist;\\
$\vdots$\\
(k) there exists a unique continuous function
$\l_k=\l_k(\l_{k+1},\cdots,\l_m,\mu,a_1,\cdots,a_k)$, which defined
on
$\{(\l_{k+1},\cdots,\l_m,\mu,a_1,\cdots,a_k):\l_{k+1}\in\mathbb{R},\cdots,\l_m\in\mathbb{R},
\mu>0,a_1\in\mathbb{R},\cdots,a_k\in\mathbb{R},a_1^2+\cdots+a_k^2<1\}$,
such that
\begin{equation*}
R_k(r,\l,\mu)\lt|_{
\begin{subarray}{l}
\l_1=\l_1(\l_2,\cdots,\l_m,\mu,a_1)\\
\l_2=\l_2(\l_3,\cdots,\l_m,\mu,a_1,a_2)\\
\cdots\cdots\\
\l_k=\l_k(\l_{k+1},\cdots,\l_m,\mu,a_1,\cdots,a_k)
\end{subarray}
}\rt.=a_k,
\end{equation*}
and
$\frac{\p\l_k(\l_{k+1},\cdots,\l_m,\mu,a_1,\cdots,a_k)}{\p\l_{k+1}},
\cdots,\frac{\p\l_k(\l_{k+1},\cdots,\l_m,\mu,a_1,\cdots,a_k)}{\p\l_m},
\frac{\p\l_k(\l_{k+1},\cdots,\l_m,\mu,a_1,\cdots,a_k)}{\p\mu}$
exist,\\
then\\
(1) if $k\leq m-2$, then there exists a unique continuous function
$$\l_{k+1}=\l_{k+1}(\l_{k+2},\cdots,\l_m,\mu,a_1,\cdots,a_{k+1}),$$
which defined on
$\{(\l_{k+2},\cdots,\l_m,\mu,a_1,\cdots,a_{k+1}):\l_{k+2}\in\mathbb{R},\cdots,
\l_m\in\mathbb{R},\mu>0,a_1\in\mathbb{R},\cdots,a_{k+1}\in\mathbb{R},
a_1^2+\cdots+a_{k+1}^2<1\}$, such that
\begin{equation*}
R_{k+1}(r,\l,\mu)\lt|_{
\begin{subarray}{l}
\l_1=\l_1(\l_2,\cdots,\l_m,\mu,a_1)\\
\l_2=\l_2(\l_3,\cdots,\l_m,\mu,a_1,a_2)\\
\cdots\cdots\\
\l_k=\l_k(\l_{k+1},\cdots,\l_m,\mu,a_1,\cdots,a_k)\\
\l_{k+1}=\l_{k+1}(\l_{k+2},\cdots,\l_m,\mu,a_1,\cdots,a_{k+1})
\end{subarray}
}\rt.=a_{k+1},
\end{equation*}
and
$\frac{\p\l_{k+1}(\l_{k+2},\cdots,\l_m,\mu,a_1,\cdots,a_{k+1})}{\p\l_{k+2}},
\cdots,\frac{\p\l_{k+1}(\l_{k+2},\cdots,\l_m,\mu,a_1,\cdots,a_{k+1})}{\p\l_m},
\frac{\p\l_{k+1}(\l_{k+2},\cdots,\l_m,\mu,a_1,\cdots,a_{k+1})}{\p\mu}$
exist;\\
(2) if $k=m-1$, then there exists a unique continuous function
$\l_m=\l_m(\mu,a_1,\cdots,a_m)$, which defined on
$\{(\mu,a_1,\cdots,a_m):\mu>0,a_1\in\mathbb{R},\cdots,a_m\in\mathbb{R},
a_1^2+\cdots+a_m^2<1\}$, such that
\begin{equation*}
R_{k+1}(r,\l,\mu)\lt|_{
\begin{subarray}{l}
\l_1=\l_1(\l_2,\cdots,\l_m,\mu,a_1)\\
\l_2=\l_2(\l_3,\cdots,\l_m,\mu,a_1,a_2)\\
\cdots\cdots\\
\l_{m-1}=\l_{m-1}(\l_m,\mu,a_1,\cdots,a_{m-1})\\
\l_m=\l_m(\mu,a_1,\cdots,a_m)
\end{subarray}
}\rt.=a_m.
\end{equation*}

Now we will prove the result above. For $1\leq k\leq m-1$, let
\begin{equation*}
\l^*=(\l^*_1,\l^*_2,\cdots,\l^*_k,\l_{k+1},\cdots,\l_m)=\l\lt|_{
\begin{subarray}{l}
\l_1=\l_1(\l_2,\cdots,\l_m,\mu,a_1)\\
\l_2=\l_2(\l_3,\cdots,\l_m,\mu,a_1,a_2)\\
\cdots\cdots\\
\l_k=\l_k(\l_{k+1},\cdots,\l_m,\mu,a_1,\cdots,a_k)
\end{subarray}
}\rt.,
\end{equation*}
where
\begin{equation*}
\begin{split}
&\l_1^*=\l_1\lt|_{
\begin{subarray}{l}
\l_1=\l_1(\l_2,\cdots,\l_m,\mu,a_1)\\
\l_2=\l_2(\l_3,\cdots,\l_m,\mu,a_1,a_2)\\
\cdots\cdots\\
\l_k=\l_k(\l_{k+1},\cdots,\l_m,\mu,a_1,\cdots,a_k)
\end{subarray}
}\rt., \l_2^*=\l_2\lt|_{
\begin{subarray}{l}
\l_2=\l_2(\l_3,\cdots,\l_m,\mu,a_1,a_2)\\
\cdots\cdots\\
\l_k=\l_k(\l_{k+1},\cdots,\l_m,\mu,a_1,\cdots,a_k)
\end{subarray}
}\rt.,\\
&\cdots, \l_k^*=\l_k\lt|_{
\begin{subarray}{l}
\l_k=\l_k(\l_{k+1},\cdots,\l_m,\mu,a_1,\cdots,a_k)
\end{subarray}
}\rt..
\end{split}
\end{equation*}
Consider the function $R_{k+1}(r,\l^*,\mu)$. A simple calculation
gives
\begin{equation}\label{Rk+1p}
\frac{\p R_{k+1}(r,\l^*,\mu)}{\p\l_{k+1}}=\lt.\lt(R_{(k+1)1}\frac{\p
\l_1^*}{\p\l_{k+1}}+R_{(k+1)2}\frac{\p \l_2^*}{\p\l_{k+1}}+\cdots
+R_{(k+1)k}\frac{\p
\l_k^*}{\p\l_{k+1}}+R_{(k+1)(k+1)}\rt)\rt|_{\l=\l^*}.
\end{equation}
By the condition $(1)-(k)$, we have for $j=1,2,\cdots,k$,
$R_j(r,\l^*,\mu)=a_j$ and consequently $\frac{\p
R_j(r,\l^*,\mu)}{\p\l_{k+1}}=0$, which is
\begin{equation}\label{Rjp}
\lt.\lt(R_{j1}\frac{\p \l_1^*}{\p\l_{k+1}}+R_{j2}\frac{\p
\l_2^*}{\p\l_{k+1}}+\cdots +R_{jk}\frac{\p
\l_k^*}{\p\l_{k+1}}+R_{j(k+1)}\rt)\rt|_{\l=\l^*}=0\ \ \mbox{for
$j=1,2,\cdots,k$}.
\end{equation}
By \eqref{Rjp}, \eqref{Rk+1p} and Lemma \ref{lemmafcz}, we have
\begin{equation*}
\frac{\p R_{k+1}(r,\l^*,\mu)}{\p\l_{k+1}}=\lt.\frac{\lt|
\begin{array}{cccc}
R_{11}&R_{12}&\cdots&R_{1(k+1)}\\
R_{21}&R_{22}&\cdots&R_{2(k+1)}\\
\multicolumn{4}{c}\dotfill\\
R_{(k+1)1}&R_{(k+1)2}&\cdots&R_{(k+1)(k+1)}\end{array}\rt|}{\lt|
\begin{array}{cccc}
R_{11}&R_{12}&\cdots&R_{1k}\\
R_{21}&R_{22}&\cdots&R_{2k}\\
\multicolumn{4}{c}\dotfill\\
R_{k1}&R_{k2}&\cdots&R_{kk}\end{array}\rt|}\rt|_{\l=\l^*}.
\end{equation*}
Then by \eqref{1}, we obtain $\frac{\p
R_{k+1}(r,\l^*,\mu)}{\p\l_{k+1}}<0$, which shows that
$R_{k+1}(r,\l^*,\mu)$ is strictly decreasing as a function of
$\l_{k+1}$. Note that $-1<R_{k+1}(r,\l^*,\mu)<1$ and
$R_{k+1}(r,\l^*,\mu)$ is bounded since (i) and (ii) in Step 1. Thus,
$R_{k+1}(r,\l^*,\mu)$, as a function of $\l_{k+1}$, respectively has
finite limit as $\l_{k+1}\to+\iy$ and as $\l_{k+1}\to-\iy$.

We claim that $R_{k+1}(r,\l^*,\mu)\to-\sqrt{1-a_1^2-\cdots-a_k^2}$
as $\l_{k+1}\to+\iy$, and
$R_{k+1}(r,\l^*,\mu)\to\sqrt{1-a_1^2-\cdots-a_k^2}$ as
$\l_{k+1}\to-\iy$.

As $\l_{k+1}\to+\iy$. Note that for $j=1,2,\cdots,k+1$,
$\lt.\frac{-\frac{\l_j}{\l_{k+1}}}{\sqrt{\sum_{i=1}^{k+1}\lt(\frac{\l_i}{\l_{k+1}}\rt)^2}}
\rt|_{\l=\l^*}$ is bounded since
$\lt|\lt.\frac{-\frac{\l_j}{\l_{k+1}}}{\sqrt{\sum_{i=1}^{k+1}\lt(\frac{\l_i}{\l_{k+1}}\rt)^2}}
\rt|_{\l=\l^*}\rt|\leq1$. Then there exists a subsequence
$\lt(\l_{k+1}\rt)_p\to+\iy$ such that for $j=1,2,\cdots,k+1$,
$\lt.\frac{-\frac{\l_j}{\l_{k+1}}}{\sqrt{\sum_{i=1}^{k+1}\lt(\frac{\l_i}{\l_{k+1}}\rt)^2}}
\rt|_{\l=\l^*|_{\l_{k+1}=\lt(\l_{k+1}\rt)_p}}$ has a finite limit
$t_j$. Let $\lt(\l^*\rt)_p=\l^*|_{\l_{k+1}=\lt(\l_{k+1}\rt)_p}$.
Then we have
\begin{equation}\label{jx1}
\lim_{p\to\iy}\lt.\frac{-\frac{\l_j}{\l_{k+1}}}
{\sqrt{\sum_{i=1}^{k+1}\lt(\frac{\l_i}{\l_{k+1}}\rt)^2}}
\rt|_{\l=\lt(\l^*\rt)_p}=t_j\ \ \mbox{for $j=1,2,\cdots,k+1$}.
\end{equation}
We only need to prove that
$R_{k+1}(r,\lt(\l^*\rt)_p,\mu)\to-\sqrt{1-a_1^2-\cdots-a_k^2}$ as
$p\to\iy$. Let
$(A(\o))_p=((A_1(\o))_p,\cdots,(A_m(\o))_p)=A_{r,\lt(\l^*\rt)_p,\mu}(\o)$.
By \eqref{Ady} and \eqref{jx1} we obtain for $j=1,2,\cdots,k+1$,
\begin{equation}\label{jx2}
\begin{split}
\lim_{p\to\iy}\frac{(A_j(\o))_p}{\sqrt{1+|(A(\o))_p|^2}}
&=\lim_{p\to\iy}\lt.\frac{\frac1\mu\lt(\frac{1}{|rN-\o|^n}l_j-\l_j\rt)}
{\sqrt{1+\frac1{\mu^2}\sum^m_{i=1}\lt(\frac{1}{|rN-\o|^n}l_i-\l_i\rt)^2}}\rt|
_{\l=\lt(\l^*\rt)_p}\\
&=\lim_{p\to\iy}\lt.\frac{\frac{1}{|rN-\o|^n}\frac{l_j}{\l_{k+1}}-\frac{\l_j}{\l_{k+1}}}
{\sqrt{\frac{\mu^2}{\l^2_{k+1}}+\sum^m_{i=1}\lt(\frac{1}{|rN-\o|^n}
\frac{l_i}{\l_{k+1}}-\frac{\l_i}{\l_{k+1}}\rt)^2}}\rt|
_{\l=\lt(\l^*\rt)_p}\\
&=\lim_{p\to\iy}\lt.\frac{-\frac{\l_j}{\l_{k+1}}}{\sqrt{\sum_{i=1}^{k+1}\lt(\frac{\l_i}{\l_{k+1}}\rt)^2}}
\rt|_{\l=\lt(\l^*\rt)_p}=t_j
\end{split}
\end{equation}
uniformly for $\o\in S$. By the Lebesgue's dominated convergence
theorem and \eqref{RIdy}, \eqref{jx2}, we have for
$j=1,2,\cdots,k+1$,
\begin{equation}\label{jx3}
\begin{split}
\lim_{p\to\iy}R_j(r,\lt(\l^*\rt)_p,\mu)
=\lim_{p\to\iy}\int_S\frac{(A_j(\o))_p}{\sqrt{1+|(A(\o))_p|^2}}\,d\s
=\int_S\lim_{p\to\iy}\frac{(A_j(\o))_p}{\sqrt{1+|(A(\o))_p|^2}}\,d\s
=t_j.
\end{split}
\end{equation}
Note that $R_j(r,\lt(\l^*\rt)_p,\mu)\equiv a_j$ for $j=1,2,\cdots,k$
by the condition (1)-(k), and $\sum_{j=1}^{k+1}t_j^2=1$,
$t_{k+1}\leq0$ by \eqref{jx2}. Then by \eqref{jx3} we have $t_j=a_j$
for $j=1,2,\cdots,k$, and $t_{k+1}=-\sqrt{1-a_1^2-\cdots-a_k^2}$.
Consequently
$$\lim_{p\to\iy}R_{k+1}(r,\lt(\l^*\rt)_p,\mu)=-\sqrt{1-a_1^2-\cdots-a_k^2}.$$
The first claim is proved.

Using the method of the proof of the first claim, we can prove the
second claim. It is proved that $R_{k+1}(r,\l^*,\mu)$ is continuous
and strictly decreasing from $\sqrt{1-a_1^2-\cdots-a_k^2}$ to
$-\sqrt{1-a_1^2-\cdots-a_k^2}$ as $\l_{k+1}$ increasing from $-\iy$
to $+\iy$. Thus, for any
$$-\sqrt{1-a_1^2-\cdots-a_k^2}<a_{k+1}<\sqrt{1-a_1^2-\cdots-a_k^2}$$
and $a_1^2+\cdots+a_k^2<1$, we have that\\
(1) if $k\leq m-2$, then there exists a unique real number
$\l_{k+1}(\l_{k+2},\cdots,\l_m,\mu,a_1,\cdots,a_{k+1})$ such that
$$R_{k+1}(r,\l^*,\mu)|_{\l_{k+1}=
\l_{k+1}(\l_{k+2},\cdots,\l_m,\mu,a_1,\cdots,a_{k+1})}=a_{k+1};$$
further, using the implicit function theorem, we have that the
function\\
$\l_{k+1}(\l_{k+2},\cdots,\l_m,\mu,a_1,\cdots,a_{k+1})$ defined on
on
$\{(\l_{k+2},\cdots,\l_m,\mu,a_1,\cdots,a_{k+1}):\l_{k+2}\in\mathbb{R},\cdots,$
$\l_m\in\mathbb{R},\mu>0,a_1\in\mathbb{R},\cdots,a_{k+1}\in\mathbb{R},
a_1^2+\cdots+a_{k+1}^2<1\}$ is a continuous function, and\\
$\frac{\p\l_{k+1}(\l_{k+2},\cdots,\l_m,\mu,a_1,\cdots,a_{k+1})}{\p\l_{k+2}},
\cdots,\frac{\p\l_{k+1}(\l_{k+2},\cdots,\l_m,\mu,a_1,\cdots,a_{k+1})}{\p\l_m},
\frac{\p\l_{k+1}(\l_{k+2},\cdots,\l_m,\mu,a_1,\cdots,a_{k+1})}{\p\mu}$
exist;\\
(2) if $k=m-1$, then there exists a unique real number
$\l_m(\mu,a_1,\cdots,a_m)$ such that
$$R_{k+1}(r,\l^*,\mu)|_{\l_m=
\l_m(\mu,a_1,\cdots,a_m)}=a_m;$$ further, using the implicit
function theorem, we have that the function
$\l_m(\mu,a_1,\cdots,a_m)$ defined on on
$\{(\mu,a_1,\cdots,a_m):\mu>0,a_1\in\mathbb{R},\cdots,a_m\in\mathbb{R},
a_1^2+\cdots+a_m^2<1\}$ is a continuous function.

\underline{Step 5}: For the case that $m\geq2$, by Step 3 and Step 4
we have that there exists a unique continuous mapping
\begin{equation*}
\l(\mu,a)=\l\lt|_{
\begin{subarray}{l}
\l_1=\l_1(\l_2,\cdots,\l_m,\mu,a_1)\\
\cdots\cdots\\
\l_k=\l_k(\l_{k+1},\cdots,\l_m,\mu,a_1,\cdots,a_k)\\
\cdots\cdots\\
\l_m=\l_m(\mu,a_1,\cdots,a_m)
\end{subarray}
}\rt.
\end{equation*}
defined on
$\{(\mu,a):\mu>0,a\in\mathbb{R}^m,a=(a_1,\cdots,a_m),|a|^2<1\}$,
such that
$$
R_{j}(r,\l(\mu,a),\mu)=a_j\ \ \ \mbox{for $j=1,2,\cdots,m$},
$$
and $\frac{\p\l_1(\mu,a)}{\p\mu},\cdots,\frac{\p\l_m(\mu,a)}{\p\mu}$
exist, where $(\l_1(\mu,a),\cdots,\l_m(\mu,a))=\l(\mu,a)$.

\underline{Step 6}: For $m\geq1$, by Step 2 and Step 5 we know that
there exists a unique continuous mapping $\l(\mu,a)$ defined on
$\{(\mu,a):\mu>0,a\in\mathbb{R}^m,|a|^2<1\}$, such that
\begin{equation}\label{Rdya1}
R(r,\l(\mu,a),\mu)=a,
\end{equation}
and $\frac{\p\l_1(\mu,a)}{\p\mu},\cdots,\frac{\p\l_m(\mu,a)}{\p\mu}$
exist, where $(\l_1(\mu,a),\cdots,\l_m(\mu,a))=\l(\mu,a)$.

In the following, we consider the function $I(r,\l_(\mu,a),\mu)$.

For a fixed $a=(a_1,\cdots,a_m)\in\mathbb{R}^m$ with $|a|^2<1$,
write $$\l(\mu,a)=\l(\mu)=(\l_1(\mu),\cdots,\l_m(\mu)).$$ Then
\begin{equation}\label{Ids}
\frac{dI(r,\l(\mu),\mu)}{d\mu}
=\lt.\lt(I_1\frac{d\l_1(\mu)}{d\mu}+I_2\frac{d\l_2(\mu)}{d\mu}
+\cdots+I_m\frac{d\l_m(\mu)}{d\mu}+I_\mu\rt)\rt|_{\l=\l(\mu)}.
\end{equation}
By \eqref{Rdya1}, we know that
\begin{equation}\label{Rdya2}
R_{j}(r,\l(\mu),\mu)=a_j\ \ \ \mbox{for $j=1,2,\cdots,m$}
\end{equation}
and
\begin{equation}\label{Rdyads}
\lt.\lt(R_{j1}\frac{d\l_1(\mu)}{d\mu}+R_{j2}\frac{d\l_2(\mu)}{d\mu}
+\cdots+R_{jm}\frac{d\l_m(\mu)}{d\mu}
+R_{j\mu}\rt)\rt|_{\l=\l(\mu)}=0\ \ \ \mbox{for $j=1,2,\cdots,m$}.
\end{equation}
Then by \eqref{Rdyads}, \eqref{Ids} and Lemma \ref{lemmafcz}, we
have
\begin{equation*}
\frac{dI(r,\l(\mu),\mu)}{d\mu}=\lt.\frac{\lt|
\begin{array}{ccccc}
R_{11}&R_{12}&\cdots&R_{1m}&R_{1\mu}\\
R_{21}&R_{22}&\cdots&R_{2m}&R_{2\mu}\\
\multicolumn{4}{c}\dotfill\\
R_{m1}&R_{m2}&\cdots&R_{mm}&R_{m\mu}\\
I_1&I_2&\cdots&I_m&I_\mu\end{array}\rt|}{\lt|
\begin{array}{cccc}
R_{11}&R_{12}&\cdots&R_{1m}\\
R_{21}&R_{22}&\cdots&R_{2m}\\
\multicolumn{4}{c}\dotfill\\
R_{m1}&R_{m2}&\cdots&R_{mm}\end{array}\rt|}\rt|_{\l=\l(\mu)}.
\end{equation*}
By \eqref{2}, we have $\frac{dI(r,\l(\mu),\mu)}{d\mu}>0$, which
shows that $I(r,\l(\mu),\mu)$ is strictly increasing as a function
of $\mu$. By (iii) in Step 1, we know that $I(r,\l(\mu),\mu)$
respectively has finite limit as $\mu\to0$ and as $\mu\to+\iy$.

We claim that $I(r,\l(\mu),\mu)\to0$ as $\mu\to0$, and
$I(r,\l(\mu),\mu)\to\sqrt{1-|a|^2}$ as $\mu\to+\iy$.

As $\mu\to0$, there exists a subsequence $\mu_k\to0$ such that
$\l_1(\mu_k)$ has a finite limit $t$ or tend to $\iy$. We only need
to prove that $I(r,\l(\mu_k),\mu_k)\to0$ as $k\to\iy$. Since
$I(r,\l(\mu_k),\mu_k)=\int_S\frac{1}
{\sqrt{1+|A_{r,\l(\mu_k),\mu_k}(\o)|^2}}\,d\s$, we only need to
prove that $|A_{r,\l(\mu_k),\mu_k}(\o)|\to+\iy$ almost everywhere on
$S$. Note that
$$|A_{r,\l(\mu_k),\mu_k}(\o)|=
\frac{1}{\mu_k}\lt|\frac{1}{|rN-\o|^n}l-\l(\mu_k)\rt|
\geq\frac{1}{\mu_k}\lt|\frac{1}{|rN-\o|^n}-\l_1(\mu_k)\rt|$$ and
$$\frac{1}{(1+r)^n}\leq\frac{1}{|rN-\o|^n}\leq\frac{1}{(1-r)^n}.$$
If $\l_1(\mu_k)\to t$ as $k\to\iy$, then
$\frac{1}{|rN-\o|^n}-\l_1(\mu_k)$ is bounded and
$\frac{1}{|rN-\o|^n}-\l_1(\mu_k)\neq0$ almost everywhere on $S$.
Thus $|A_{r,\l(\mu_k),\mu_k}(\o)|\to+\iy$ almost everywhere on $S$.
If $\l_1(\mu_k)\to\iy$ as $k\to\iy$, then it is obvious that
$|A_{r,\l(\mu_k),\mu_k}(\o)|\to+\iy$ uniformly for $\o\in S$. The
first claim is proved.

As $\mu\to+\iy$, $\frac{1}{\mu}\frac{1}{|rN-\o|^n}\to0$ uniformly
for $\o\in S$. For $j=1$ or $j=2$ or $\cdots$ or $j=m$ , if there
exists a subsequence $\mu_k\to+\iy$ such that
$\l_j(\mu_k)/\mu_k\to\iy$, then $|A_{r,\l(\mu_k),\mu_k}(\o)|\to+\iy$
uniformly for $\o\in S$, and $I(r,\l(\mu_k),\mu_k)\to 0$, a
contradiction. This shows that for $j=1,2,\cdots,m$, $\l_j(\mu)/\mu$
is bounded as $\mu\to+\iy$. Thus there exists a subsequence
$\mu_k\to+\iy$ such that $-\l_j(\mu_k)/\mu_k$ tend to a finite limit
$t_j$ for $j=1,2,\cdots,m$. That is
\begin{equation}\label{lmu}
\lim_{k\to\iy}-\l_j(\mu_k)/\mu_k=t_j\ \ \ \mbox{for
$j=1,2,\cdots,m$}.
\end{equation}
we only need to prove that $I(r,\l(\mu_k),\mu_k)\to\sqrt{1-|a|^2}$
as $k\to\iy$. Let
$$(A(\o))_k=((A_1(\o))_k,\cdots,(A_m(\o))_k)=A_{r,\l(\mu_k),\mu_k}(\o).$$
By \eqref{Ady} and \eqref{lmu} we obtain for $j=1,2,\cdots,m$,
\begin{equation}\label{wR}
\begin{split}
\lim_{k\to\iy}\frac{(A_j(\o))_k}{\sqrt{1+|(A(\o))_k|^2}}
&=\lim_{k\to\iy}\frac{\frac1\mu_k\lt(\frac{1}{|rN-\o|^n}l_j-\l_j(\mu_k)\rt)}
{\sqrt{1+\frac1{\mu_k^2}\sum^m_{i=1}\lt(\frac{1}{|rN-\o|^n}l_i-\l_i(\mu_k)\rt)^2}}\\
&=\lim_{k\to\iy}\frac{-\frac{\l_j(\mu_k)}{\mu_k}}
{\sqrt{1+\sum^m_{i=1}\lt(\frac{\l_i(\mu_k)}{\mu_k}\rt)^2}}
=\frac{t_j}{\sqrt{1+\sum^m_{i=1}t_i^2}}
\end{split}
\end{equation}
uniformly for $\o\in S$, and
\begin{equation}\label{wI}
\begin{split}
\lim_{k\to\iy}\frac{1}{\sqrt{1+|(A(\o))_k|^2}}
&=\lim_{k\to\iy}\frac{1}
{\sqrt{1+\frac1{\mu_k^2}\sum^m_{i=1}\lt(\frac{1}{|rN-\o|^n}l_i-\l_i(\mu_k)\rt)^2}}\\
&=\lim_{k\to\iy}\frac{1}
{\sqrt{1+\sum^m_{i=1}\lt(\frac{\l_i(\mu_k)}{\mu_k}\rt)^2}}
=\frac{1}{\sqrt{1+\sum^m_{i=1}t_i^2}}
\end{split}
\end{equation}
uniformly for $\o\in S$. By the Lebesgue's dominated convergence
theorem and \eqref{RIdy}, \eqref{wR}, \eqref{wI} we have for
$j=1,2,\cdots,m$,
\begin{equation}\label{wjR}
\begin{split}
\lim_{k\to\iy}R_j(r,\l(\mu_k),\mu_k)
&=\lim_{k\to\iy}\int_S\frac{(A_j(\o))_k}{\sqrt{1+|(A(\o))_k|^2}}\,d\s\\
&=\int_S\lim_{k\to\iy}\frac{(A_j(\o))_k}{\sqrt{1+|(A(\o))_k|^2}}\,d\s\\
&=\frac{t_j}{\sqrt{1+\sum^m_{i=1}t_i^2}},
\end{split}
\end{equation}
and
\begin{equation}\label{wjI}
\begin{split}
\lim_{k\to\iy}I(r,\l(\mu_k),\mu_k)
&=\lim_{k\to\iy}\int_S\frac{1}{\sqrt{1+|(A(\o))_k|^2}}\,d\s\\
&=\int_S\lim_{k\to\iy}\frac{1}{\sqrt{1+|(A(\o))_k|^2}}\,d\s\\
&=\frac{1}{\sqrt{1+\sum^m_{i=1}t_i^2}}.
\end{split}
\end{equation}
Note that $R_j(r,\l(\mu_k),\mu_k)\equiv a_j$ for $j=1,2,\cdots,m$ by
\eqref{Rdya2}, and
$$\sum^m_{j=1}\lt(\frac{t_j}{\sqrt{1+\sum^m_{i=1}t_i^2}}\rt)^2
+\lt(\frac{1}{\sqrt{1+\sum^m_{i=1}t_i^2}}\rt)^2=1.$$ Then by
\eqref{wjR} we obtain that
$\frac{t_j}{\sqrt{1+\sum^m_{i=1}t_i^2}}=a_j$ for $j=1,2,\cdots,m$,
and $\frac{1}{\sqrt{1+\sum^m_{i=1}t_i^2}}=\sqrt{1-|a|^2}$.
Consequently by \eqref{wjI},
$$\lim_{k\to\iy}I(r,\l(\mu_k),\mu_k)=\sqrt{1-|a|^2}.$$
The second claim is proved.

It is proved that $I(r,\l(\mu),\mu)$ is continuous and strictly
increasing from $0$ to $\sqrt{1-|a|^2}$ as $\mu$ increasing from 0
to $+\iy$. Thus, for any $0<b<\sqrt{1-|a|^2}$ and $|a|<1$, there
exists a unique real number $\mu(a,b)$ such that
$I(r,\l(\mu(a,b)),\mu(a,b))=b$. Further, using the implicit function
theorem, we have the function $\mu(a,b)$ defined on
$\{(a,b):a\in\mathbb{R}^m,b\in\mathbb{R},|a|<1,\
0<b<\sqrt{1-|a|^2}\}$ is a continuous function.

Denote $\l(\mu(a,b))$ by $\l(r,a,b)$. Denote $\mu(a,b)$ by
$\mu(r,a,b)$. We have proved that there exist a unique pair of
continuous mappings $\l=\l(r,a,b)$ and $\mu=\mu(r,a,b)$ such that
$R(r,\l(r,a,b),\mu(r,a,b))=a$ and $I(r,\l(r,a,b),\mu(r,a,b))=b$ on
the upper half ball. The lemma is proved. \qed

Now we give the proof of Lemma 2.

\noindent{\it Proof of Lemma 2.}\quad We will prove Lemma 2 by two
cases: $m=1$ and $m\geq2$. The case that $m=1$ will be proved in
Step 1. The case that $m\geq2$ will be proved in Step 2 - Step 5.

\underline{Step 1}: For the case that $m=1$, we have
\begin{equation*}
\mc R(r,\l)=\int_S\frac{\frac{1}{|rN-\o|^n}-\l}
{|\frac{1}{|rN-\o|^n}-\l|}\,d\s=
\begin{cases}
1, &\l\in(-\iy,\frac{1}{(1+r)^n}],\\
\int_S\frac{\frac{1}{|rN-\o|^n}-\l}
{|\frac{1}{|rN-\o|^n}-\l|}\,d\s, &\l\in(\frac{1}{(1+r)^n},\frac{1}{(1-r)^n}),\\
-1, &\l\in[\frac{1}{(1-r)^n},+\iy).
\end{cases}
\end{equation*}
Obviously $\mc R(r,\l)\equiv1$ when $\l\leq\frac{1}{(1+r)^n}$, $\mc
R(r,\l)\equiv-1$ when $\l\geq\frac{1}{(1-r)^n}$, and $\mc R(r,\l)$
is continuous and strictly decreasing from $1$ to $-1$ as $\l$
increasing from $\frac{1}{(1+r)^n}$ to $\frac{1}{(1-r)^n}$. Then for
any $-1<a<1$, there exists a unique real number $\l(a)$ such that
\begin{equation*}
\mc R(r,\l)\lt|_{\l=\l(a)}\rt.=a.
\end{equation*}
Further, using the implicit function theorem, we have that the
function $\l=\l(a)$ defined on $\{a: -1<a<1\}$ is a continuous
function. Write $\l(a)=\l(r,a)$. Then the case that $m=1$ is proved.

\underline{Step 2}: For the case that $m\geq2$, we give some
denotation and calculation. Let
$$\mc A_{r,\l}(\omega)=\mc A(\o)=(\mc A_1(\o),\mc A_2(\o),\cdots,\mc A_m(\o)),$$
$$\mc R(r,\l)=(\mc R_1(r,\l),\mc R_2(r,\l),\cdots,\mc R_m(r,\l)),$$
$$l=(l_1,\cdots,l_m), \l=(\l_1,\l_2,\cdots,\l_m), \ \mbox{and} \ \
a=(a_1,a_2,\cdots,a_m).$$ By \eqref{A}, $|\mc
A(\o)|=\sqrt{\lt(\frac{1}{|rN-\o|^n}-\l_1\rt)^2+\l_2^2+\cdots+\l_m^2}$.
So if let set
$$H=\{\l=(\l_1,\cdots,\l_m)\in\mathbb{R}^m: \l_2=\cdots=\l_m=0\},$$
then obviously for $i,j=1,2,\cdots,m$, $\frac{\p \mc
R_j(r,\l)}{\p\l_i}$ exist for $\l\in\mathbb{R}^m\backslash H$. We
denote $\frac{\p \mc R_j(r,\l)}{\p\l_i}=\mc R_{ji}$ for
$i,j=1,2,\cdots,m$. Then a simple calculation gives that for
$\l\in\mathbb{R}^m\backslash H$,
\begin{equation}\label{Rjj}
\mc R_{jj}=-\int_S \frac{|\mc A(\o)|^2-A^2_j(\o)}{|\mc A(\o)|^3}d\s\
\ \ \ \mbox{for $j=1,2,\cdots,m$};
\end{equation}
\begin{equation}\label{Rji}
\mc R_{ji}=-\int_S \frac{-\mc A_i(\o)\mc A_j(\o)}{|\mc A(\o)|^3}d\s\
\ \ \ \mbox{for $i\neq j; i,j=1,2,\cdots,m$};
\end{equation}

It is easy to see that\\
(1) by \eqref{A} and \eqref{R}, for $j=1,2,\cdots,m$, $\mc
R_j(r,\l)$ is a continuous function for any $\l\in\mathbb{R}^m$;\\
(2) by \eqref{A} and \eqref{R}, for $j=1,2,\cdots,m$, fixing the
components of $\l$ expect $\l_j$, $\mc R_j(r,\l)\to-1$ or $1$
according to $\l_j\to+\iy$ or $\l_j\to-\iy$;\\
(3) by \eqref{A} and \eqref{Rjj}, $\mc R_{11}<0$ for any
$\l\in\mathbb{R}^m\backslash H$, and $\mc R_1(r,\l)$ is strictly
decreasing as a function of $\l_1$ for fixed $\l_2,\cdots,\l_m$ with
$\l_2,\cdots,\l_m$ are not all $0$;\\
(4) by \eqref{R}, for fixed $\l_2=\cdots=\l_m=0$,
\begin{equation*}
\mc R_1(r,\l)=\int_S\frac{\frac{1}{|rN-\o|^n}-\l_1}
{|\frac{1}{|rN-\o|^n}-\l_1|}\,d\s=
\begin{cases}
1, &\l_1\in(-\iy,\frac{1}{(1+r)^n}],\\
\int_S\frac{\frac{1}{|rN-\o|^n}-\l_1}
{|\frac{1}{|rN-\o|^n}-\l_1|}\,d\s, &\l_1\in(\frac{1}{(1+r)^n},\frac{1}{(1-r)^n}),\\
-1, &\l_1\in[\frac{1}{(1-r)^n},+\iy);
\end{cases}
\end{equation*}\\
(5) by \eqref{A} and \eqref{Rjj}, for $j=2,\cdots,m$, $\mc R_{jj}<0$
for any $\l\in\mathbb{R}^m$, and $\mc R_j(r,\l)$ is strictly
decreasing as a function of $\l_j$ for fixed the other components of
$\l$.

In addition, for $\l\in\mathbb{R}^m\backslash H$, we claim that
\begin{equation}\label{h1}
\frac{\lt|
\begin{array}{cccc}
\mc R_{11}&\mc R_{12}&\cdots&\mc R_{1(k+1)}\\
\mc R_{21}&\mc R_{22}&\cdots&\mc R_{2(k+1)}\\
\multicolumn{4}{c}\dotfill\\
\mc R_{(k+1)1}&\mc R_{(k+1)2}&\cdots&\mc
R_{(k+1)(k+1)}\end{array}\rt|}{\lt|
\begin{array}{cccc}
\mc R_{11}&\mc R_{12}&\cdots&\mc R_{1k}\\
\mc R_{21}&\mc R_{22}&\cdots&\mc R_{2k}\\
\multicolumn{4}{c}\dotfill\\
\mc R_{k1}&\mc R_{k2}&\cdots&\mc R_{kk}\end{array}\rt|}<0\ \ \
\mbox{for integer $k$ with $1\leq k\leq m-1$}.
\end{equation}
Now we will prove the claim above.

For \eqref{Rjj} and \eqref{Rji}, let $d\tilde{\s}=(1/|\mc
A(\o)|^3)d\s$, $T=\int_S d\tilde{\s}$, $d\xi=(1/T)d\tilde{\s}$,
$\tilde{b}=\int_S |\mc A(\o)|^2d\xi$, and for $i,j=1,2,\cdots,m$,
$\tilde{a}_{ij}=\int_S -\mc A_i(\o)\mc A_j(\o)d\xi$, $c_j=\int_S \mc
A_j(\o)d\xi$. Then $T>0$, $\int_S d\xi=1$, and
\begin{equation}\label{hRjjdy2}
\mc R_{jj}=-\frac{T}{\mu}(\tilde{b}+\tilde{a}_{jj})\ \ \ \ \mbox{for
$j=1,2,\cdots,m$};
\end{equation}
\begin{equation}\label{hRjidy2}
\mc R_{ji}=-\frac{T}{\mu}\tilde{a}_{ij}\ \ \ \ \mbox{for $i\neq j;
i,j=1,2,\cdots,m$}.
\end{equation}
Since $\mc A_1(\o)=\frac{1}{|rN-\o|^n}-\l_1$ by \eqref{A} and
$\int_S d\xi=1$, we have
\begin{equation}\label{ha11c1}
\begin{split}
-\tilde{a}_{11}-c_1^2&=\int_S A_1^2(\o)d\xi-\lt(\int_S A_1(\o)d\xi\rt)^2\\
&=\int_S \lt[A_1(\o)-\int_S A_1(\o)d\xi\rt]^2d\xi>0.
\end{split}
\end{equation}
Since $\int_S d\xi=1$ and $\mc A_j(\o)=-\l_j$ for $j=2,\cdots,m$ by
\eqref{A}, we have
\begin{equation}\label{haij}
\tilde{a}_{ij}=-\int_S \mc A_i(\o)d\xi\int_S \mc
A_j(\o)d\xi=-c_ic_j\ \ \ \ \mbox{for $i\neq1$ or $j\neq1$},
i,j=1,2,\cdots,m;
\end{equation}
and
\begin{equation}\label{hxhls}
\lt|\begin{array}{cc}
\tilde{a}_{11}&\tilde{a}_{1j}\\
\tilde{a}_{j1}&\tilde{a}_{jj}
\end{array}\rt|=\lt|\begin{array}{cc}
\tilde{a}_{11}&-c_1c_j\\
-c_jc_1&-c_jc_j
\end{array}\rt|
=c_j^2(-\tilde{a}_{11}-c_1^2)=\l_j^2(-\tilde{a}_{11}-c_1^2)\mbox{
for $j=2,\cdots,m$}.
\end{equation}

For integer $1\leq p\leq m$, let
\begin{equation}\label{hQp}
Q_p=\lt|
\begin{array}{cccc}
\tilde{b}+\tilde{a}_{11}&\tilde{a}_{12}&\cdots&\tilde{a}_{1p}\\
\tilde{a}_{21}&\tilde{b}+\tilde{a}_{22}&\cdots&\tilde{a}_{2p}\\
\vdots&\vdots&\ddots&\vdots\\
\tilde{a}_{p1}&\tilde{a}_{p2}&\cdots&\tilde{b}+\tilde{a}_{pp}\end{array}
\rt|.
\end{equation}
Since $\l\in\mathbb{R}^m\backslash H$, we have that when $p=1$,
$Q_1=\tilde{b}+\tilde{a}_{11}>0$. By \eqref{haij} and Lemma
\ref{lemmaQn}, we have that when $p\geq2$,
\begin{equation*}
\begin{split}
Q_p&=\tilde{b}^p+\tilde{b}^{p-1}\sum_{j=1}^p\tilde{a}_{jj}
+\tilde{b}^{p-2}\sum_{j=2}^p\lt|
\begin{array}{cc}
\tilde{a}_{11}&\tilde{a}_{1j}\\
\tilde{a}_{j1}&\tilde{a}_{jj}\end{array}\rt|\\
&=\tilde{b}^{p-1}\lt(\tilde{b}+\sum_{j=1}^p\tilde{a}_{jj}\rt)
+\tilde{b}^{p-2}\sum_{j=2}^p \lt|
\begin{array}{cc}
\tilde{a}_{11}&\tilde{a}_{1j}\\
\tilde{a}_{j1}&\tilde{a}_{jj}\end{array}\rt|.
\end{split}
\end{equation*}
Consequently, by \eqref{hxhls}, \eqref{ha11c1}, $\tilde{b}>0$, and
$\l\in\mathbb{R}^m\backslash H$, we obtain that when $p\geq2$,
$Q_p>0$.

By \eqref{hRjjdy2}, \eqref{hRjidy2} and \eqref{hQp}, we have for
integer $k$ with $1\leq k\leq m-1$,
\begin{equation*}
\begin{split}
\frac{\lt|
\begin{array}{cccc}
\mc R_{11}&\mc R_{12}&\cdots&\mc R_{1(k+1)}\\
\mc R_{21}&\mc R_{22}&\cdots&\mc R_{2(k+1)}\\
\multicolumn{4}{c}\dotfill\\
\mc R_{(k+1)1}&\mc R_{(k+1)2}&\cdots&\mc
R_{(k+1)(k+1)}\end{array}\rt|}{\lt|
\begin{array}{cccc}
\mc R_{11}&\mc R_{12}&\cdots&\mc R_{1k}\\
\mc R_{21}&\mc R_{22}&\cdots&\mc R_{2k}\\
\multicolumn{4}{c}\dotfill\\
\mc R_{k1}&\mc R_{k2}&\cdots&\mc R_{kk}\end{array}\rt|}
&=\frac{\lt(-T\rt)^{k+1} \lt|
\begin{array}{cccc}
\tilde{b}+\tilde{a}_{11}&\tilde{a}_{12}&\cdots&\tilde{a}_{1(k+1)}\\
\tilde{a}_{21}&\tilde{b}+\tilde{a}_{22}&\cdots&\tilde{a}_{2(k+1)}\\
\vdots&\vdots&\ddots&\vdots\\
\tilde{a}_{(k+1)1}&\tilde{a}_{(k+1)2}&\cdots&\tilde{b}+\tilde{a}_{(k+1)(k+1)}\end{array}
\rt|}{\lt(-T\rt)^k\lt|
\begin{array}{cccc}
\tilde{b}+\tilde{a}_{11}&\tilde{a}_{12}&\cdots&\tilde{a}_{1k}\\
\tilde{a}_{21}&\tilde{b}+\tilde{a}_{22}&\cdots&\tilde{a}_{2k}\\
\vdots&\vdots&\ddots&\vdots\\
\tilde{a}_{k1}&\tilde{a}_{k2}&\cdots&\tilde{b}+\tilde{a}_{kk}\end{array}
\rt|}\\
&=\frac{\lt(-T\rt)^{k+1}Q_{k+1}}{\lt(-T\rt)^{k}Q_k}\\
&=\lt(-T\rt)\frac{Q_{k+1}}{Q_k}.
\end{split}
\end{equation*}
Note that $T>0, Q_k>0, Q_{k+1}>0$. Then the claim \eqref{h1} is
proved.

\underline{Step 3}: For the case that $m=2$, by (1)-(3) in Step 2,
we know that for fixed $\l_2,\cdots,\l_m$ with $\l_2,\cdots,\l_m$
are not all $0$, $\mc R_1(r,\l)$ is strictly decreasing from $1$ to
$-1$ as $\l_1$ increasing from $-\iy$ to $+\iy$. by (4) in Step 2,
we know that for fixed $\l_2=\cdots=\l_m=0$, $\mc R_1(r,\l)\equiv1$
when $\l_1\leq\frac{1}{(1+r)^n}$, $\mc R_1(r,\l)\equiv-1$ when
$\l_1\geq\frac{1}{(1-r)^n}$, and $\mc R_1(r,\l)$ is continuous and
strictly decreasing from $1$ to $-1$ as $\l_1$ increasing from
$\frac{1}{(1+r)^n}$ to $\frac{1}{(1-r)^n}$. Then for any $-1<a_1<1$
and any fixed $\l_2,\cdots,\l_m$, there exists a unique real number
$\l_1(\l_2,\cdots,\l_m,a_1)$ such that
\begin{equation*}
\mc R_1(r,\l)\lt|_{\l_1=\l_1(\l_2,\cdots,\l_m,a_1)}\rt.=a_1.
\end{equation*}
Further, using the implicit function theorem, we have that the
function $\l_1=\l_1(\l_2,\cdots,\l_m,a_1)$ defined on
$\{(\l_2,\cdots,\l_m,a_1):\l_2\in\mathbb{R},\cdots,\l_m\in\mathbb{R},-1<a_1<1\}$
is a continuous function and
$\frac{\p\l_1(\l_2,\cdots,\l_m,a_1)}{\p\l_2},
\cdots,\frac{\p\l_1(\l_2,\cdots,\l_m,a_1)}{\p\l_m}$ exist for
$(\l_2,\cdots,\l_m,a_1)$ with $\l_2,\cdots,\l_m$ are not all $0$.

\underline{Step 4}: For the case that $m\geq2$, we will prove the
following result:

For an integer $k$ with $1\leq k\leq m-1$, if\\
(1) there exists a unique continuous function
$\l_1=\l_1(\l_2,\cdots,\l_m,a_1)$, which defined on\\
$\{(\l_2,\cdots,\l_m,a_1):\l_2\in\mathbb{R},\cdots,\l_m\in\mathbb{R},-1<a_1<1\}$,
such that
\begin{equation*}
\mc R_1(r,\l)\lt|_{\l_1=\l_1(\l_2,\cdots,\l_m,a_1)}\rt.=a_1,
\end{equation*}
and $\frac{\p\l_1(\l_2,\cdots,\l_m,a_1)}{\p\l_2},
\cdots,\frac{\p\l_1(\l_2,\cdots,\l_m,a_1)}{\p\l_m}$ exist for
$(\l_2,\cdots,\l_m,a_1)$ with $\l_2,\cdots,\l_m$ are not all $0$;\\
(2) there exists a unique continuous function
$\l_2=\l_2(\l_3,\cdots,\l_m,a_1,a_2)$, which defined on\\
$\{(\l_3,\cdots,\l_m,a_1,a_2):\l_3\in\mathbb{R},\cdots,\l_m\in\mathbb{R},
a_1\in\mathbb{R},a_2\in\mathbb{R},a_1^2+a_2^2<1\}$, such that
\begin{equation*}
\mc R_2(r,\l)\lt|_{\l_1=\l_1(\l_2,\cdots,\l_m,a_1)
\atop\l_2=\l_2(\l_3,\cdots,\l_m,a_1,a_2)}\rt.=a_2,
\end{equation*}
and $\frac{\p\l_2(\l_3,\cdots,\l_m,a_1,a_2)}{\p\l_3},
\cdots,\frac{\p\l_2(\l_3,\cdots,\l_m,a_1,a_2)}{\p\l_m}$ exist for
$(\l_3,\cdots,\l_m,a_1,a_2)$ with $\l_3,\cdots,\l_m$ are not all
$0$;\\
$\vdots$\\
(k) there exists a unique continuous function
$\l_k=\l_k(\l_{k+1},\cdots,\l_m,a_1,\cdots,a_k)$, which defined on
$\{(\l_{k+1},\cdots,\l_m,a_1,\cdots,a_k):\l_{k+1}\in\mathbb{R},\cdots,\l_m\in\mathbb{R},
a_1\in\mathbb{R},\cdots,a_k\in\mathbb{R},a_1^2+\cdots+a_k^2<1\}$,
such that
\begin{equation*}
\mc R_k(r,\l)\lt|_{
\begin{subarray}{l}
\l_1=\l_1(\l_2,\cdots,\l_m,a_1)\\
\l_2=\l_2(\l_3,\cdots,\l_m,a_1,a_2)\\
\cdots\cdots\\
\l_k=\l_k(\l_{k+1},\cdots,\l_m,a_1,\cdots,a_k)
\end{subarray}
}\rt.=a_k,
\end{equation*}
and $\frac{\p\l_k(\l_{k+1},\cdots,\l_m,a_1,\cdots,a_k)}{\p\l_{k+1}},
\cdots,\frac{\p\l_k(\l_{k+1},\cdots,\l_m,a_1,\cdots,a_k)}{\p\l_m}$
exist for $(\l_{k+1},\cdots,\l_m,a_1,\cdots,a_k)$ with\\
$\l_{k+1},\cdots,\l_m$ are not all $0$,\\
then\\
(1) if $k\leq m-2$, then there exists a unique continuous function\\
$\l_{k+1}=\l_{k+1}(\l_{k+2},\cdots,\l_m,a_1,\cdots,a_{k+1})$, which
defined on
$\{(\l_{k+2},\cdots,\l_m,a_1,\cdots,a_{k+1}):\l_{k+2}\in\mathbb{R},\cdots,
\l_m\in\mathbb{R},a_1\in\mathbb{R},\cdots,a_{k+1}\in\mathbb{R},a_1^2+\cdots+a_{k+1}^2<1\}$,
such that
\begin{equation*}
\mc R_{k+1}(r,\l)\lt|_{
\begin{subarray}{l}
\l_1=\l_1(\l_2,\cdots,\l_m,a_1)\\
\l_2=\l_2(\l_3,\cdots,\l_m,a_1,a_2)\\
\cdots\cdots\\
\l_k=\l_k(\l_{k+1},\cdots,\l_m,a_1,\cdots,a_k)\\
\l_{k+1}=\l_{k+1}(\l_{k+2},\cdots,\l_m,a_1,\cdots,a_{k+1})
\end{subarray}
}\rt.=a_{k+1},
\end{equation*}
and
$\frac{\p\l_{k+1}(\l_{k+2},\cdots,\l_m,a_1,\cdots,a_{k+1})}{\p\l_{k+2}},
\cdots,\frac{\p\l_{k+1}(\l_{k+2},\cdots,\l_m,a_1,\cdots,a_{k+1})}{\p\l_m}$
exist for\\
$(\l_{k+2},\cdots,\l_m,a_1,\cdots,a_{k+1})$ with
$\l_{k+2},\cdots,\l_m$ are not all $0$;\\
(2) if $k=m-1$, then there exists a unique continuous function
$\l_{m}=\l_{m}(a_1,\cdots,a_m)$, which defined on
$\{(a_1,\cdots,a_m):a_1\in\mathbb{R},\cdots,a_m\in\mathbb{R},a_1^2+\cdots+a_m^2<1\}$,
such that
\begin{equation*}
\mc R_m(r,\l)\lt|_{
\begin{subarray}{l}
\l_1=\l_1(\l_2,\cdots,\l_m,a_1)\\
\l_2=\l_2(\l_3,\cdots,\l_m,a_1,a_2)\\
\cdots\cdots\\
\l_{m-1}=\l_{m-1}(\l_m,\cdots,\l_m,a_1,\cdots,a_{m-1})\\
\l_m=\l_m(a_1,\cdots,a_m)
\end{subarray}
}\rt.=a_m.
\end{equation*}

Now we will prove the result above. For $1\leq k\leq m-1$, let
\begin{equation*}
\l^*=(\l^*_1,\l^*_2,\cdots,\l^*_k,\l_{k+1},\cdots,\l_m)=\l\lt|_{
\begin{subarray}{l}
\l_1=\l_1(\l_2,\cdots,\l_m,a_1)\\
\l_2=\l_2(\l_3,\cdots,\l_m,a_1,a_2)\\
\cdots\cdots\\
\l_k=\l_k(\l_{k+1},\cdots,\l_m,a_1,\cdots,a_k)
\end{subarray}
}\rt.,
\end{equation*}
where
\begin{equation*}
\begin{split}
&\l_1^*=\l_1\lt|_{
\begin{subarray}{l}
\l_1=\l_1(\l_2,\cdots,\l_m,a_1)\\
\l_2=\l_2(\l_3,\cdots,\l_m,a_1,a_2)\\
\cdots\cdots\\
\l_k=\l_k(\l_{k+1},\cdots,\l_m,a_1,\cdots,a_k)
\end{subarray}
}\rt., \l_2^*=\l_2\lt|_{
\begin{subarray}{l}
\l_2=\l_2(\l_3,\cdots,\l_m,a_1,a_2)\\
\cdots\cdots\\
\l_k=\l_k(\l_{k+1},\cdots,\l_m,a_1,\cdots,a_k)
\end{subarray}
}\rt.,\\
&\cdots, \l_k^*=\l_k\lt|_{
\begin{subarray}{l}
\l_k=\l_k(\l_{k+1},\cdots,\l_m,a_1,\cdots,a_k)
\end{subarray}
}\rt..
\end{split}
\end{equation*}
Consider the function $\mc R_{k+1}(r,\l^*)$. A simple calculation
gives for $\l^*$ with $\l_{k+1},\cdots,\l_m$ are not all $0$,
\begin{equation}\label{HRk+1p}
\frac{\p \mc R_{k+1}(r,\l^*)}{\p\l_{k+1}}=\lt.\lt(\mc
R_{(k+1)1}\frac{\p \l_1^*}{\p\l_{k+1}}+\mc R_{(k+1)2}\frac{\p
\l_2^*}{\p\l_{k+1}}+\cdots +\mc R_{(k+1)k}\frac{\p
\l_k^*}{\p\l_{k+1}}+\mc R_{(k+1)(k+1)}\rt)\rt|_{\l=\l^*}.
\end{equation}
By the condition (1)-(k), we have for $j=1,2,\cdots,k$, $\mc
R_j(r,\l^*)=a_j$ and consequently $\frac{\p \mc
R_j(r,\l^*)}{\p\l_{k+1}}=0$ for $\l^*$ with $\l_{k+1},\cdots,\l_m$
are not all $0$, which is
\begin{equation}\label{HRjp}
\lt.\lt(\mc R_{j1}\frac{\p \l_1^*}{\p\l_{k+1}}+\mc R_{j2}\frac{\p
\l_2^*}{\p\l_{k+1}}+\cdots +\mc R_{jk}\frac{\p
\l_k^*}{\p\l_{k+1}}+\mc R_{j(k+1)}\rt)\rt|_{\l=\l^*}=0\ \ \mbox{for
$j=1,2,\cdots,k$}.
\end{equation}
By \eqref{HRjp}, \eqref{HRk+1p} and Lemma \ref{lemmafcz}, we have
\begin{equation*}
\frac{\p \mc R_{k+1}(r,\l^*)}{\p\l_{k+1}}=\lt.\frac{\lt|
\begin{array}{cccc}
\mc R_{11}&\mc R_{12}&\cdots&\mc R_{1(k+1)}\\
\mc R_{21}&\mc R_{22}&\cdots&\mc R_{2(k+1)}\\
\multicolumn{4}{c}\dotfill\\
\mc R_{(k+1)1}&\mc R_{(k+1)2}&\cdots&\mc
R_{(k+1)(k+1)}\end{array}\rt|}{\lt|
\begin{array}{cccc}
\mc R_{11}&\mc R_{12}&\cdots&\mc R_{1k}\\
\mc R_{21}&\mc R_{22}&\cdots&\mc R_{2k}\\
\multicolumn{4}{c}\dotfill\\
\mc R_{k1}&\mc R_{k2}&\cdots&\mc
R_{kk}\end{array}\rt|}\rt|_{\l=\l^*}
\end{equation*}
for $\l^*$ with $\l_{k+1},\cdots,\l_m$ are not all $0$. Then by
\eqref{h1}, we obtain $\frac{\p \mc R_{k+1}(r,\l^*)}{\p\l_{k+1}}<0$
for $\l^*$ with $\l_{k+1},\cdots,\l_m$ are not all $0$, which shows
that when $\l_{k+1}\neq0$, $\mc R_{k+1}(r,\l^*)$ is strictly
decreasing as a function of $\l_{k+1}$. Since $\mc R_{k+1}(r,\l^*)$
is continuous as a function of $\l_{k+1}$ by the condition (1)-(k)
and (1) in Step 2, then for $\l_{k+1}\in\mathbb{R}$, $\mc
R_{k+1}(r,\l^*)$ is strictly decreasing as a function of $\l_{k+1}$.
Note that $-1<\mc R_{k+1}(r,\l^*)<1$ and $\mc R_{k+1}(r,\l^*)$ is
bounded since (2) and (5) in Step 2. Thus $\mc R_{k+1}(r,\l^*)$, as
a function of $\l_{k+1}$, respectively has finite limit as
$\l_{k+1}\to+\iy$ and as $\l_{k+1}\to-\iy$.

We claim that $\mc R_{k+1}(r,\l^*)\to-\sqrt{1-a_1^2-\cdots-a_k^2}$
as $\l_{k+1}\to+\iy$, and\\ $\mc
R_{k+1}(r,\l^*)\to\sqrt{1-a_1^2-\cdots-a_k^2}$ as $\l_{k+1}\to-\iy$.

As $\l_{k+1}\to+\iy$. Note that for $j=1,2,\cdots,k+1$,
$\lt.\frac{-\frac{\l_j}{\l_{k+1}}}{\sqrt{\sum_{i=1}^{k+1}\lt(\frac{\l_i}{\l_{k+1}}\rt)^2}}
\rt|_{\l=\l^*}$ is bounded since
$\lt|\lt.\frac{-\frac{\l_j}{\l_{k+1}}}{\sqrt{\sum_{i=1}^{k+1}\lt(\frac{\l_i}{\l_{k+1}}\rt)^2}}
\rt|_{\l=\l^*}\rt|\leq1$. Then there exists a subsequence
$\lt(\l_{k+1}\rt)_p\to+\iy$ such that for $j=1,2,\cdots,k+1$,
$\lt.\frac{-\frac{\l_j}{\l_{k+1}}}{\sqrt{\sum_{i=1}^{k+1}\lt(\frac{\l_i}{\l_{k+1}}\rt)^2}}
\rt|_{\l=\l^*|_{\l_{k+1}=\lt(\l_{k+1}\rt)_p}}$ has a finite limit
$t_j$. Let $\lt(\l^*\rt)_p=\l^*|_{\l_{k+1}=\lt(\l_{k+1}\rt)_p}$.
Then we have
\begin{equation}\label{hjx1}
\lim_{p\to\iy}\lt.\frac{-\frac{\l_j}{\l_{k+1}}}
{\sqrt{\sum_{i=1}^{k+1}\lt(\frac{\l_i}{\l_{k+1}}\rt)^2}}
\rt|_{\l=\lt(\l^*\rt)_p}=t_j\ \ \mbox{for $j=1,2,\cdots,k+1$}.
\end{equation}
We only need to prove that $\mc
R_{k+1}(r,\lt(\l^*\rt)_p)\to-\sqrt{1-a_1^2-\cdots-a_k^2}$ as
$p\to\iy$. Let $$(\mc A(\o))_p=((\mc A_1(\o))_p,\cdots,(\mc
A_m(\o))_p)=\mc A_{r,\lt(\l^*\rt)_p}(\o).$$ By \eqref{A} and
\eqref{hjx1} we obtain for $j=1,2,\cdots,k+1$,
\begin{equation}\label{hjx2}
\begin{split}
\lim_{p\to\iy}\frac{(\mc A_j(\o))_p}{|(\mc A(\o))_p|}
&=\lim_{p\to\iy}\lt.\frac{\frac{1}{|rN-\o|^n}l_j-\l_j}
{\sqrt{\sum^m_{i=1}\lt(\frac{1}{|rN-\o|^n}l_i-\l_i\rt)^2}}\rt|
_{\l=\lt(\l^*\rt)_p}\\
&=\lim_{p\to\iy}\lt.\frac{\frac{1}{|rN-\o|^n}\frac{l_j}{\l_{k+1}}-\frac{\l_j}{\l_{k+1}}}
{\sqrt{\sum^m_{i=1}\lt(\frac{1}{|rN-\o|^n}
\frac{l_i}{\l_{k+1}}-\frac{\l_i}{\l_{k+1}}\rt)^2}}\rt|
_{\l=\lt(\l^*\rt)_p}\\
&=\lim_{p\to\iy}\lt.\frac{-\frac{\l_j}{\l_{k+1}}}{\sqrt{\sum_{i=1}^{k+1}\lt(\frac{\l_i}{\l_{k+1}}\rt)^2}}
\rt|_{\l=\lt(\l^*\rt)_p}=t_j
\end{split}
\end{equation}
uniformly for $\o\in S$. By the Lebesgue's dominated convergence
theorem and \eqref{R}, \eqref{hjx2}, we have for $j=1,2,\cdots,k+1$,
\begin{equation}\label{hjx3}
\begin{split}
\lim_{p\to\iy}\mc R_j(r,\lt(\l^*\rt)_p)
=\lim_{p\to\iy}\int_S\frac{(\mc A_j(\o))_p}{|(\mc A(\o))_p|}\,d\s
=\int_S\lim_{p\to\iy}\frac{(\mc A_j(\o))_p}{|(\mc A(\o))_p|}\,d\s
=t_j.
\end{split}
\end{equation}
Note that $\mc R_j(r,\lt(\l^*\rt)_p)\equiv a_j$ for $j=1,2,\cdots,k$
by the condition, and $\sum_{j=1}^{k+1}t_j^2=1$, $t_{k+1}\leq0$ by
\eqref{hjx2}. Then by \eqref{hjx3} we have $t_j=a_j$ for
$j=1,2,\cdots,k$, and $t_{k+1}=-\sqrt{1-a_1^2-\cdots-a_k^2}$.
Consequently
$$\lim_{p\to\iy}\mc R_{k+1}(r,\lt(\l^*\rt)_p)=-\sqrt{1-a_1^2-\cdots-a_k^2}.$$
The first claim is proved.

Using the method of the proof of the first claim, we can prove the
second claim. It is proved that $\mc R_{k+1}(r,\l^*)$ is continuous
and strictly decreasing from $\sqrt{1-a_1^2-\cdots-a_k^2}$ to
$-\sqrt{1-a_1^2-\cdots-a_k^2}$ as $\l_{k+1}$ increasing from $-\iy$
to $+\iy$. Therefore, for any
$$-\sqrt{1-a_1^2-\cdots-a_k^2}<a_{k+1}<\sqrt{1-a_1^2-\cdots-a_k^2}$$
with $a_1^2-\cdots-a_k^2<1$, we obtain \\
(1) if $k\leq m-2$, then there exists a unique real number
$\l_{k+1}(\l_{k+2},\cdots,\l_m,a_1,\cdots,a_{k+1})$ such that
$$\mc R_{k+1}(r,\l^*)|_{\l_{k+1}=
\l_{k+1}(\l_{k+2},\cdots,\l_m,a_1,\cdots,a_{k+1})}=a_{k+1};$$
further, using the implicit function theorem, we have that the
function $$\l_{k+1}(\l_{k+2},\cdots,\l_m,a_1,\cdots,a_{k+1})$$
defined on on
$\{(\l_{k+2},\cdots,\l_m,a_1,\cdots,a_{k+1}):\l_{k+2}\in\mathbb{R},\cdots,
\l_m\in\mathbb{R}>0,a_1\in\mathbb{R},\cdots,a_{k+1}\in\mathbb{R},
a_1^2+\cdots+a_{k+1}^2<1\}$ is a continuous function, and
$$\frac{\p\l_{k+1}(\l_{k+2},\cdots,\l_m,a_1,\cdots,a_{k+1})}{\p\l_{k+2}},
\cdots,\frac{\p\l_{k+1}(\l_{k+2},\cdots,\l_m,a_1,\cdots,a_{k+1})}{\p\l_m}$$
exist;\\
(2) if $k=m-1$, then there exists a unique real number
$\l_{m}=\l_{m}(a_1,\cdots,a_m)$ such that
$$\mc R_m(r,\l^*)|_{\l_m=
\l_{m}(a_1,\cdots,a_m)}=a_m;$$ further, using the implicit function
theorem, we have that the function\\
$\l_{m}=\l_{m}(a_1,\cdots,a_m)$ defined on
$\{(a_1,\cdots,a_m):a_1\in\mathbb{R},\cdots,a_m\in\mathbb{R},
a_1^2+\cdots+a_m^2<1\}$ is a continuous function.

\underline{Step 5}: For the case that $m\geq2$, by Step 3 and Step 4
we have that there exists a unique continuous mapping
\begin{equation*}
\l(a)=\l\lt|_{
\begin{subarray}{l}
\l_1=\l_1(\l_2,\cdots,\l_m,a_1)\\
\cdots\cdots\\
\l_k=\l_k(\l_{k+1},\cdots,\l_m,a_1,\cdots,a_k)\\
\cdots\cdots\\
\l_m=\l_m(a_1,\cdots,a_m)
\end{subarray}
}\rt.
\end{equation*}
defined on $\{a=(a_1,\cdots,a_m)\in\mathbb{R}^m:|a|<1\}$, such that
$$
R(r,\l(a))=a.
$$
Write $\l(a)=\l(r,a)$. Then the case that $m\geq2$ is proved.
 \qed

\end{document}